\documentstyle[12pt,leqno]{article}
\def\ba#1{\begin{array}{#1}}\def\ea{\end{array}}
\def\be{\begin{equation}}\def\ee{\end{equation}}
\def\ben{\begin{enumerate}}\def\en{\end{enumerate}}
\def\bi{\begin{itemize}}\def\ei{\end{itemize}}
\newtheorem{exx}{}[section]
\def\bex{\begin{exx}{\bf{Example. }}\em}\def\ex{\rm\end{exx}}
\newtheorem{sex}{}[exx]
\def\bs#1{\begin{sex}{\bf{#1}}\em}\def\es{\end{sex}\rm}
\newtheorem{ssex}{}[sex]
\def\bss{\begin{ssex}\em}\def\ess{\end{ssex}\rm}
\def\bp{\vskip15pt\begin{exx}{\bf{Proposition. }}}\def\ep{\rm\end{exx}}
\def\bl{\vskip15pt\begin{exx}{\bf{Lemma. }}}\def\el{\rm\end{exx}}
\def\schub{{\sc schubert}}
\renewcommand{\O}{\mbox{${\cal O}$}}
\renewcommand{\L}{\mbox{${\cal L}$}}
\def\M{{\cal M}}
\def\A{{\cal A}}
\def\B{{\cal B}}
\def\C{{\cal C}}
\def\J{{\cal J}}
\def\R{{\cal R}}
\def\D{{\cal D}}
\def\I{{\bf I}}\def\l{{\bf L}}
\def\Q{\mbox{${\cal Q}$}}
\def\cs{\mbox{${\cal S}$}}
\def\n{\noindent}
\def\pf{\vskip10pt\n{\bf Proof.\hskip10pt}}

\def\p#1{\mbox{I$\!$P$^{#1}$}}
\def\pp#1{\mbox{\footnotesize I$\!$P$^{#1}$}}
\def\q{\mbox{$\p1\!\times\!\p1$}}
\def\qp{\mbox{\footnotesize$\pp1\!\!\times\!\pp1$}}
\def\ar{\mbox{${\rightarrow}$}}
\def\sear#1{\mbox{$\searrow\hskip-5pt\raise4pt\hbox{$#1$}$}}
\def\lar{\mbox{${\longrightarrow}$}}
\def\m{{\underline{m}}}
\def\eml#1#2{\mbox{${{\cal E}({#1},{#2})}$} }
\def\up#1{\mbox{$^{\underline{\hbox{\footnotesize{#1}}}}$} }
\def\bin#1#2{\mbox{$(\stackrel{#1}{_{#2}})$}}
\def\qed{\hfill\mbox{$\Box$}}
\def\z {\mbox{$Z$}}
\def\inj {\mbox{\raise-.03cm\hbox{\mbox{\tiny\bf$^\succ$} }
\hskip-0.25cm$\ar$}}
\def\surj{\mbox{$-\!\!\ar$\hskip-.08cm\raise-.03cm\hbox
{\mbox{\tiny\bf$^\succ$}}}}
\def\injdown{\mbox{\hskip.01cm\raise.035cm
\hbox{\mbox{\large$\check{}$}}
\hskip-.245cm$\downarrow$}}
\def\surjdown
{\mbox{\hskip.01cm\raise-.37cm\hbox{\mbox{\large$\check{}$}}
\hskip-.248cm$\downarrow$}}
\begin{document}
{\Large\bf
\centerline{Enumeration of $n$-fold tangent}
\centerline{ hyperplanes to a surface}
}
\centerline{\bf by}
\centerline{\bf Israel Vainsencher
\footnote{Partially supported by CNPq, Brasil.}}
\centerline{\bf Universidade Federal de Pernambuco}
{\small
\centerline{ Cidade Universit\'aria, 50670--901, Recife -- Pe, Brasil}
\centerline{israel@dmat.ufpe.br}
}
\begin{abstract} {\small\n For each $1\leq n\leq6$ we present
formulas for the number of $n-$nodal curves in an $n-$dimensional
linear system on a smooth, projective surface. This yields in
particular the numbers of rational curves in the system of hyperplane
sections of a generic $K3-$surface imbedded in \p{n} by a complete
system of curves of genus $n$ as well as the number
{\bf17,601,000} of rational ({\em singular}) plane quintic
curves in a generic quintic threefold. } \end{abstract}
\section{Introduction} \normalsize
The purpose of this article is to present formulas for the number of
$n-$nodal curves in an $n-$dimensional linear system on a smooth,
projective surface for  $1\leq n\leq6$.
The method also yields formulas for the number of multi--tangent planes
to a hypersurface. In particular, it enables us to find the number
{\bf17,601,000}  of rational ({\em singular}) plane {quintic} curves
in a generic {quintic} threefold. We give several examples and
discuss the difficulties involved  for $n\geq7$.

Our motivation was in response to a question asked by A. Lopez and C.
Ciliberto regarding the number of rational curves in the system of
hyperplane sections in a generic $K3-$surface imbedded in \p{4}
(resp. \p{5}) as a $(2,3)$ (resp. $(2,2,2)$)--complete intersection.
In \cite{clm} (joint with Miranda) they study degenerations $K3\ar$
union of  2 scrolls. According to A. Lopez (priv. comm.), the
consideration of limit curves in the scrolls suggests a formula for
the number of rational curves in the $K3-$surface. However, the
numbers they have found are so far in disagreement with
those obtained by the formulas presented here for $n=4,5,6$ (cf. \ref{k3}).

A similar question communicated by S. Katz concerns the number of
{\em plane} rational curves of degree 5 contained in a generic
quintic $3-$fold in \p{4}.

The subject was raised by Clemens \cite{c} and has received
striking contributions from physicists (cf. Morrison \cite{mor},
Piene \cite{ragni}, Bershadski {\it et al.} \cite{vafa}). The total number of
rational curves of degree $\geq4$ has not been verified so far
{\it au go\^ut du jour des math\'ematiciens}. The cases of
degrees  1 and  2 were treated by  Harris \cite{har} and Katz
\cite{katz}. As for degree  3, it required a thorough
investigation of the Chow ring of the variety of twisted cubics
(cf. Ellingsrud and Str{\o}mme \cite{twc},\cite{twc2} (see also
\cite{twc3} for a simpler approach). A pleasant byproduct was the
development of the computer package \schub\ by Katz and Str{\o}mme
\cite{schub}.

The work of Coray \cite{co} reduces certain enumerative questions
concerning rational curves in \p3 to the question of finding the
numbers $\Delta_{\mu,\nu}$ of irreducible rational curves of
bidegree $(\mu,\nu)$ passing through $2(\mu+\nu)-1$ general points
on a quadric surface.
He computes $\Delta_{2,3}$ and $\Delta_{2,4}$ (in addition to a
few trivial cases).
We also obtain here $\Delta_{3,3}$ (\ref{d33}),
$\Delta_{2,5}$ (\ref{d25}), and $\Delta_{3,4}$ (\ref{d34}).

 Counting  hyperplanes multi--tangent to a curve is well known as a
particular case of the classical formula of De Jonqui\`eres \cite{harris},
\cite{i}.

For surfaces, the cases $n\leq 3$ are classical and have been checked with
currently standard tools of intersection theory, cf. Kleiman
\cite{kleimansing},\cite{i}. The degrees of the ``Severi
varieties'' of nodal curves in the plane were computed
(in principle) by Ran in \cite{ranb},\cite{ranb1}.

Although we have at our disposal multiple point formulas (Kleiman
\cite{kleimanmultpts}, Ran \cite{ran}), they do not give the correct
answer for multi--tangencies already for $n=$2 or 3 due to the presence
of cusps. There are also
formulas taking into account stationary multiple-points (Colley
\cite{susan}). However, for $n\geq4$ the relevant map does {\em not}
satisfy a required curvilinearity hypothesis. This is due to the
existence of curves with a triple point in virtually any linear
system of dimension $\geq4$ on a surface.

Our approach is based on the iteration procedure presented in
\cite{i0},\cite{i} (also explored in a broader context in
\cite{kleimanmultpts},\cite{ran},\cite{susan}).

We obtain, for each $n=1,\dots,6$ a formula for the degree of a
zero cycle supported on the set of sequences $(C,y_1,\dots,y_n)$
such that $C$ is a member of a (sufficiently general) linear
system of dimension $n$ and $y_1$ is a singular point of $C$,
$y_{2}$ is a singular point of the blowup of $C$ at $y_1$, and
so on (roughly speaking, cf. \S2 and (\ref{step}),(\ref{dejonq})
for the precise statement).

The main novelty here is, essentially, detecting the contribution to that
zero cycle due to singularities worse than nodes (cf.\ref{formulas}). We
also sharpen the scope of validity of the formulas, now requiring only that
the relevant loci be finite (\ref{step}).

Thanks are due to the MSRI for the stimulating environment and to P. Aluffi,
E. Arrondo, S. Katz, A. Lopez and M. Pedreira for many pleasant conversations
and to C. Schoen for the comments following Example 4.6. I'm also indebted to
S.L. Kleiman for reading a preliminary version and helping to clarify the
proof of the Lemma \ref{reduced}. We also thank \schub \cite{schub},
for patiently allowing us to verify many examples.

\section{Notation and basic definitions}

We recall, for the reader's benefit, some definitions from \cite{i}.
Let $Y$ be a smooth variety. For each sequence of integers $\m =(m_1,
\dots,$ $ m_r)$ we say an effective divisor $D$ has a singularity of
({\em weak}) type $\m $ if the following holds:
\bi
\item there is a point $y_1$ of multiplicity$\geq m_1$ in $D$; next
\item blowup $Y$ at $y_1$, let $E_1$ denote the exceptional divisor and let
$D_1$ denote the {\em total} transform of $D$; then
\item
require that the effective divisor $D_1-m_1E_1$ have a point $y_2$
of multiplicity $\geq m_2$, and so on.
\ei
The sequence $(y_1,y_2,\dots)$ thus constructed is called a {\em
singularity of type} $\m$~ of $D$. We further say the type is {\em strict} if
all inequalities are equalities and each $y_i$ lies off the exceptional
divisor.

One may also consider $nested$ sequences $(\dots, m_i$ $(m_{i+1},$ $\dots),
\dots)$ and say a singularity is of such type if $y_i$ is of multiplicity$
\geq m_i$ and $y_{i+1}$ is infinitely near to $y_i,\ i.e.,$ lies on the
exceptional divisor besides being of multiplicity $\geq m_{i+1}$, etc. We
write $ m^{[k]}$ to indicate $k$ repetitions of $m$.

\bex
\label{triple}
Let $Y$ be a surface and $y_1$ a triple point on the curve $C$. Then of course
$C$ has a singularity of strict type $(3)$. However, if the  3 tangents are
distinct, $C$ also has a singularity of weak type $(2^{[4]})$ due to the
intersections of the strict transform of $C$ and the exceptional line $E_1$:
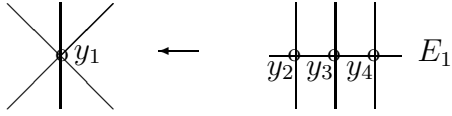
\ex
\vskip5pt
$$\begin{picture}(0,0)(0,0)
\put(-54,-4){\footnotesize o}\put(-47,-3){$y_1$}\put(-52,-2){\line( 0,-1){20}}
\put(-52,-2){\line( 0,1){20}}\put(-52,-2){\line( 1,1){20}} \put(-52,-2)
{\line( -1,-1){20}}
\put(-52,-2){\line( 1,1){20}}\put(-52,-2){\line( 1,-1){20}} \put(-52,-2)
{\line( -1,1){20}}
\put(0,0){\vector(-1,0){15}} \put(52,-2){\line( 1,0){25}}\put(52,-2)
{\line(-1,0){25}}
\put(37,-2){\line(0,1){20}}\put(37,-2){\line(0,-1){20}} \put(52,-2)
{\line( 0,1){20}}
\put(52,-2){\line(0,-1){20}}\put(67,-2){\line(0,1){20}}\put(67,-2)
{\line(0,-1){20}}
\put(34,-4){{\footnotesize o}}\put(26,-9){$y_2$}\put(49,-4){{\footnotesize o}}
\put(41,-9){$y_3$}\put(64,-4){{\footnotesize o}}\put(56,-9){$y_4$}
\put(83,-4){$E_1$}
\end{picture}$$
\vskip10pt\bex
\label{tritwo}
On the other hand, if $y_1$ is of type $(3(2))$, it follows that
$C$ has a singularity of type $(2^{[6]})$!\ex\vskip2pt
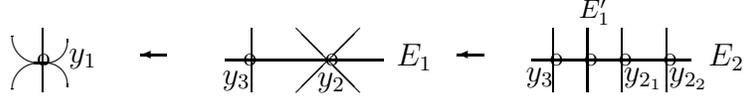
$$\begin{picture}(40,0)(0,0)
\put(-105,6){\oval(25,18)[br]}
\put(-101,7){\oval(25,20)[bl]}
\put(-105,-12){\oval(25,18)[tr]}
\put(-101,-13){\oval(25,20)[tl]}
\put(-104,-4){\footnotesize o}\put(-92,-3){$y_1$}
\put(-102,-2){\line( 0,-1){12}}
\put(-102,-2){\line( 0,1){12}}
\put(-55,0){\vector(-1,0){10}}
\put(-3,-2){\line( 1,0){30}}\put(-3,-2){\line(-1,0){30}}
\put(-23,-2){\line(0,1){12}}\put(-23,-2){\line(0,-1){12}}
\put(-26,-4){{\footnotesize o}}\put(-34,-11){$y_3$}
\put(6,-2){\line( 1,1){12}}\put(6,-2){\line(-1,-1){12}}
\put(6,-2){\line(-1,1){12}}\put(6,-2){\line(1,-1){12}}
\put(5,-4){{\footnotesize o}}\put(2,-12){$y_2$}
\put(32,-4){$E_1$}
\put(65,0){\vector(-1,0){10}}
\put(113,-2){\line( 1,0){30}}\put(113,-2){\line(-1,0){30}}
\put(91,-2){\line(0,1){12}}\put(91,-2){\line(0,-1){12}}
\put(90,-4){{\footnotesize o}}\put(81,-10){$y_3$}
\put(104,-2){\line(0,1){12}}\put(104,-2){\line(0,-1){12}}
\put(103,-4){{\footnotesize o}}\put(101,14){{\footnotesize$E_1'$}}
\put(117,-2){\line(0,1){12}}\put(117,-2){\line(0,-1){12}}
\put(116,-4){{\footnotesize o}}\put(118,-10){$y_{2_1}$}
\put(134,-2){\line(0,1){12}}\put(134,-2){\line(0,-1){12}}
\put(133,-4){{\footnotesize o}} \put(135,-10){$y_{2_2}$}
\put(150,-4){$E_2$}
 \end{picture}$$\vskip10pt\n
Indeed, let $y_2$ be the double point infinitely near to the triple
point $y_1$, and let $C_1$ denote the total transform of $C$; then
$C_1-3E_1$ is effective and intersects $E_1$ twice at $y_2$ and once
at the (smooth) branch $y_3$. Thus, the divisor $C':=C_1-2E_1$ has
multiplicity $3$ at the point $y_2$. Blowing it up, let $C_2$ be the
total transform of $C'$; now $C_2-2E_2$ still contains the
exceptional line $E_2$ once and therefore has 4 double points: one
for the intersection of $E_2$ and the strict transform of $E_1$, two
for the branches over $y_2$ and finally one over $y_3$.

\bex\label{nonred}
Let $Y$ be a surface and $y_1$
a fourfold point on the curve $C$. Then $C_1-2E_1$ is nonreduced, hence $C$
has a singularity $(y_1,\dots, y_r)$ of type $(2^{[r]})$ for any $r$.
\ex
This ilustrates a main difficulty in our approach to enumeration
of singularities. Formulas for a given type are usually not hard
to obtain, at least in principle (cf. (\ref{dejonq}) below), but
the exact contribution of each strict type actually occurring
seems less evident. For the case we're interested in, we have the
following description of the possible singularity types.
\bp
\label{tipos}
Let $Y$ be a smooth surface; fix $n\!\in\!\{1,\dots,6\}$.
Let $D$ be an ample divisor on $Y$. Then there exists $r_0$ such
that for all $r\geq r_0$ and any sufficiently general linear
subsystem $S$ of $|rD|$ of dimension $n$, there are at most
finitely many members $C\!\in\! S$ with a singularity of type
$(2^{[n]})$. Moreover, we have the following list of possible
strict types actually occurring in type $(2^{[n]})$:\vskip-5pt
$$\ba{c}
 n\leq3\Rightarrow (2^{[n]})\ only;\\
 n=4\Rightarrow (2^{[4]})\ or\ (3);\\
 n=5\Rightarrow (2^{[5]})\ or\ (3,2)\ or\ (2,3);\\
 n=6\Rightarrow (2^{[6]})\ or\ ( 3(2))\ or\ any\ of\
(3,2,2),(2,3,2),(2,2,3). \ea$$
\ep
\vskip-10pt
\pf Set $\L=\O(D)$ and let $\M_y$\ be the ideal sheaf of a point $y\!\in\!
Y$. The members of $|D|$ with an $m-$fold point at $y$ come from
$H^0(Y,\M^m_y\otimes \L)$. Let ${Y}_n^{\hbox{o}}$ denote the
complement of the diagonals in $Y^{\times{n}}$. Given a sequence of
positive integers, $(m_1,\dots,m_n)$, replacing $\L$ by a sufficiently
high power, we may assume $H^1(Y,\M^{m_1}_{y_1}\dots \M^{m_n}_{y_n}
\otimes \L)=0$ for all $(y_1,\dots,y_n)\!\in\! {Y}_n^{\hbox{o}}$. It
follows that the set $$\{(C,y_1, \dots, y_n) \!\in\!|D|\times
{Y}_n^{\hbox{o}}\ |\ \hbox{mult}_{y_i}C\geq{m_i}\}$$ is a projective
bundle over ${Y}_n^{\hbox{o}}$ with fibre dimension
 $=\hbox{dim} |D| - \Sigma m_i(m_i+1)/2$.
Its image in $|D|$ is of codimension $\Sigma m_i(m_i+1)/2-2$.
Therefore no sufficiently general subsystem of dimension$\leq3$ (resp.
$\leq7$) has a member with a triple (resp. $4-$fold)
point. It can be easily checked that a singularity
of type $(2(2))$ ($i.e.$, a double point with another infinitely
near) (resp. $(2(2),2)$ or $(2,2(2)))$ imposes  3 (resp.  4)
independent conditions.

Let $(y_1,\dots,y_6)$ be a singularity of weak type $(2^{[6]})$
occuring in a general $\infty^6$ linear system.
As explained just above, a $4-$fold point imposes 8
conditions, so each $y_i$ is at worst a triple point. Moreover,
it can be checked that  2 triple points (infinitely near or not)
impose at least  8 conditions, thus at most one of the $y_i$ is
triple. We claim that $y_i$ cannot be a triple point unless
$i\leq3$. Indeed, the imposition of 3 double points
$(y_1,y_2,y_3)$ costs at least 3 parameters, leaving less than
the  4 required for the acquisition of an additional triple point.

A similar argument rules out other sequences of double points
(with some possibly infinitely near) different from those listed.
\qed
\section{Basic setup}  \label{setup}
Let $f:X\ar S$ be proper and smooth. Let $\L$ be an invertible $\O_X-$module
and let $D\subset X$ be the scheme of zeros of a section of $\L$. As in
\cite{i}, we construct a scheme $\Sigma(\m;D)$ whose fibre over
each $s\!\in\! S$ consists of the sequences of singularities of type $\m$~ of
the fibre $D_s$.

Set $X_0=S,\ X_1=X,\ f_{1}=f:X_1\ar X_0.$ For $r\geq1$ denote by
$$b_{r+1}:X_{r+1}\ar X_{r}\times_{f_{r}}X_r
\quad\hbox{  and  }\quad
p_{r+1,i}:X_{r}\times_{f_{r}}X_r\ar{}X_r$$
respectively the blowup of the diagonal and the projection.

Set $f_{r+1,i}=p_{r+1,i}\circ{} b_{r+1}$. We think of each $X_r$ as a
scheme over $X_{r-1}$ with structure map $f_r=f_{r,1}$.

Write $E_{1,r}$ for the exceptional divisor of $b_{r}$.

For $2\leq j<r$ set $E_{r-j+1,j}=f_{r,2}^*\cdots f_{j+1,2}^*E_{1,j}$. By
abuse, still denote by the same symbol pullbacks of $E_{r-j+1,j}$ via
compositions of the structure maps $f_3,f_4,\dots$. Notice the 2\up{nd}
index in $E_{r-j+1,j}$ indicates where the divisor first appears in the
sequence of blowups, whereas $r-j$ keeps track of the number of pullbacks
via the $f_{k,2}$.

For each sequence of nonnegative integers $\m =(m_1,\dots, m_r)$ we define
the divisor on $X_{r+1}$,
$$
\m E=m_rE_{1,r+1}+\cdots +m_2E_{r-1,3}+m_1E_{r,2}.
$$
Let $y_1\!\in\! X_1$ lie over $s\!\in\! X_0$. Notice that, by
construction, the fibre $X_{2{y_1}}$ of $f_2$ over $y_1$ is
equal to the blowup of the fibre $f_1^{-1}(s)$ at $y_1$. By the
same token, a point in $X_r$ lying over $s$ should be thought of
as a sequence $(y_1,\dots,y_r)$ of points in $f_1^{-1}(s)$ each
possibly infinitely near to a previous one. Also, the fibre of
$\m E$ over a point $(y_1,\dots,y_r)\!\in\! X_r$ is equal to
$m_rE_{y_r}+\cdots +m_1E_{y_1}$, where $E_{y_i}\subset
X_{i+1y_i}$ denotes (for $i<r$, the total transform of) the
exceptional divisor of the blowup of $X_{iy_{i-1}}$ at ${y_i}$.
We set
$$
\L(\m)=f_{r+1,2}^*\cdots f_{2,2}^*\L\otimes\O_{X_{r+1}}(-\m E).
$$
Pulling back the section of $\L$ defining $D$, we get the diagram of maps of
$\O_{X_{r+1}}-$modules,
\be\label{s}
\ba{l}    \O_{X_{r+1}}\\
\downarrow\put(20,10){\vector(3,-1){40}}
\put(40,10){$\sigma_{\m}^D$}\\
\hskip-2cm f_{r+1,2}^*\cdots f_{2,2}^*\L\hskip10pt\lar\hskip10pt
f_{r+1,2}^*\cdots f_{2,2}^*\L\otimes\O_{\m E}.
\ea
\ee
\n By construction, $\sigma_{\m}^D$ vanishes on a fibre $f_{r+1}^{-1}(y_1,
\dots,y_r)$ iff $y_1,\dots,y_r$ is a singularity of type $\m$~ of $D_{s}$,
where $s= f_1(y_1)$.

We define the $\m-${\em contact sheaf} as the $\O_{X_r}-$module,
$$
\eml{\m}{\L}=f_{r+1*}(\O_{\m E}\otimes f_{r+1,2}^*\cdots f_{2,2}^*\L).
$$
\bl
\label{sh}
Notation as above, we have:
\ben
\item \eml{\m}{\L} is a locally free \O$_{X_r}-$ module of rank $\Sigma
\ (^{dim f+m_i-1}_{\ \ \ \ dim f})$ and its formation commutes with base
change;
\item there are exact sequences,
$$0\lar \eml{m_r}{\L(\m')}\lar \eml{\m}{\L} \lar f_{r}^*\eml{\m'}{\L}\lar0,$$
where $\m'$ denotes the truncated sequence $(m_1,\dots,m_{r-1});$
\item we have $\eml{1}\L=\L$ and for $\mu\geq2$ we have an exact sequence,
$$0\lar \L\otimes{}Sym^{\mu-1}\Omega^1_{X/S}\lar \eml{\mu}{\L} \lar
\eml{\mu-1}{\L} \lar0.$$
\en
\el
\pf
The inclusion $f_{r+1,2}^*\m'E\subset \m E$ yields the exact sequence
\be\label{exs}
\ba{ccc}
0\hskip3pt\lar &\hskip-15pt\O_{\m E}(-f_{r+1,2}^*\m'E)&\hskip-15pt\lar
\hskip3pt\O_{\m E} \hskip3pt\lar\hskip3pt \O_{f_{r+1,2}^*\m'E}\hskip3pt
\lar0\\
&||&\\
&\hskip-20pt\O_{m_rE_{1,r+1}}(-f_{r+1,2}^*\m'E)&
\ea
\ee
Notice $f_{r+1,2}^*\m'E$ and $m_rE_{1,r+1}$ are $f_{r+1,1}-$flat.
Indeed, for a divisor such as $E_{2,r}:=f_{r+1,2}^*E_{1,r}$ which
intersects the blowup center $\Delta(X_r)$ properly (along
$\Delta(E_{1,r})$), the total and strict transforms are one and
the same. Thus, to show $f_{r+1,1}-$ flat\-ness of $E_{2,r}$ it
suffices to verify that each power of the ideal sheaf of $\Delta(E_{1,r})$ in
$p_{r+1,2}^*E_{1,r}$ is $p_{r+1,1}-$flat. This is a consequence
of the following.
\bl
\label{flat}
Let $p:X\ar Y$ be a smooth map of smooth varieties. Let $Z\subset X$
be a smooth subvariety of $X$ such that the restriction of $p$ induces
an isomorphism $Z\stackrel{\sim}{\ar} p(Z)$ onto a hypersurface of $Y$.
Let ${\cal I}$ denote the ideal of $Z$ in $X$. Then each power ${\cal I} ^m$
is $p-$flat.
\el\pf
We assume for simplicity dim $p$=1 (hence codim($Z,X$) $=2$). There
is a local representation of $p$ by a ring homomorphism ${\cal A\ar B}$
fitting into a commutative diagram,
$$\ba{ccc}
   \A &\lar & \B\\
   \uparrow&&\uparrow\\
   \C:=\R[u]&\lar& \D:=\R[u,v]
\ea
$$
such that the vertical maps are \'etale, $\R$ is regular, $u,v$ denote
indeterminates and the image of $u$ (resp. $u,v$) generates the ideal
of $p(Z)$ (resp. $Z$) (cf. \cite{ak0}, p. 128--130). Under these
circumstances, let $\M$ be a $\D-$module flat/$\C$. Then $\B\otimes_{\D}
\M$ is flat/$\A$. Indeed, put $\A':=\A\otimes_{\C}\D$; clearly
$\M_{\A}:=\A\otimes_{\C}\M$ is an $\A'-$module flat/$\A$. Notice $\A\ar B$
factors as $\A\ar \A' \ar \B$ and $\B$ is \'etale, hence flat/$\A'$ . Let
${\J}\subset \A$ be an ideal. We have $0\ar \J\otimes_{\A}\M_{\A}\ar \M_{\A}$
exact. Hence
$$\ba{cccc}
0\ar & \B\otimes_{\A'}\J\otimes_{\A}\M_{\A}&\ar & \B\otimes_{\A'}\M_{A}\\
    &||                      &  &||             \\
0\ar &\J\otimes_{\A}\B\otimes_{\D}\M&\ar & \B\otimes_{\D}\M
\ea
$$
is exact by flatness of $\B$/$\A'$. Apply this to the ideal
$\M=(u,v)^m\D$, which is a flat, in fact free $\C-$module with basis
$\{u^m,\dots,uv^{m-1},v^m,\dots\}$.
\hfill\mbox{$\Box$(for \ref{flat})}\vskip10pt

\n The same argument applies to all $E_{j,r-j+1}$. Since a sum of
flat divisors is flat, we've proved that $\m E$ is $f_{r+1,1}-$flat.

 \n Tensoring (\ref{exs}) with $f_{r+1,2}^*\cdots
f_{2,2}^*\L$ and pushing forward by $f_{r+1}=p_{r+1,1}b_{r+1}$,
the assertions follow by a standard base change argument (cf. \cite{i}, p.
411).\vskip1pt\hfill\mbox{$\Box$(for \ref{sh})}
\bp
\label{step}
Let $\Sigma(\m;D)\subset X_r$ be the scheme of zeros of the section
$\sigma_{\m}^D:\O_{X_r}\ar \eml{\m}{\L}$ defined in {\rm(\ref{s})}. Then:
\ben
\item $\Sigma(\m;D)$ is equal to the scheme of zeros of $
\sigma_{\m}^D$ along the fibres of $f_{r+1}$, thus parametrizing the
singularities of type $\m$~ of the fibres of $D$;
\item with notation as in Lemma \ref{sh}, setting $D'=$
$f_{r+1,2}^*(f_{r,2}^*\dots f_{2,2}^*D-\m'E)$ restricted
over $\Sigma(\m';D)$, we have
$$\Sigma(\m;D)=\Sigma((m_r);D');$$
\item  each component of $\Sigma(\m;D)$ is of codimension$\leq \rho=
\Sigma\ (^{dim f+m_i-1}_{\ \ \ \ dim f});$
\item \label{mreg}
if $\Sigma(\m;D)$ is empty or of the right codimension $\rho$ then
its class in the Chow group of $X_r$ is given by the formula,
\en \be \label{dejonq}
[\Sigma(\m;D)]=c_\rho(\eml{\m}{\L})\cap[X_r].
\ee\ep

\pf The  1\up{st} assertion follows from \cite{ak}, Prop.(2.3).
The  2\up{nd} one derives from the exact sequence in
Lemma \ref{sh}(2).
The remaining are well known facts (cf. Fulton\cite{f}).\qed

\bs{Remark.}\label{defmreg}
In practice, the formula (\ref{dejonq}) may be computed using
the exact sequences in \ref{sh}. However, it is only useful to the extent
the conditions of (\ref{step})\ref{mreg} are met; we then say $D$
is $\m-${\em regular}. We refer to \cite{i} for sufficient
conditions for $\m-$regularity.\es

\bp\label{iter}
Let $D\subset X\ar S$ be as in the beginning of \S\ref{setup}.
Set $S'=\Sigma(2,S)$. Fix $P\!\in\! D$. Assume that
\ben
\item $S$ is regular at the image of $P$;
\item the ``total space'' $D$ is smooth at $P$ and
\item the fibre of $D$ through $P$ has an ordinary double
point ({\em odp}) there.
\en
Then we have that $S'$ is smooth at $P$.
Moreover, $D':=f_{2,2}^*D_{|S'}-2E_{1,2|S'}$ is regular along the
inverse image of $P$.
\ep\pf
We assume for simplicity dim$X/S=2$ and dim$S\geq1$. The question is
local analytic. Let $A$ be a regular local ring and $\M$ its
maximal ideal, let $ h\!\in\! B=A[|x_1,x_2|]$ and set ${\cal N}=\M B +
(x_1,x_2)B$. Assume that $B/(h)$ is regular and $h=x_1x_2\ mod (x_1,
x_2)^3 +{\M}B$. Then $\bar{B}:=B/(h,h_{x_1},h_{x_2})$ is regular.
Indeed, we may write $h=t+m_1x_1+m_2x_2 $ $+x_1x_2+\cdots$, with
$t,m_1,$ $m_2\!\in\!{\M}$. Notice that, since $h\!\in\!{\cal N}-{\cal N} ^2$,
we have in fact $t\!\in\!{\M}-{\M}^2$. From $h_{x_i}=m_i + x_j + \cdots
(\{i,j\} = \{1,2\})$, it follows that $h, h_{x_1}, h_{x_2}$ are linearly
independent $mod\ {\cal N}^2$, as desired for the regularity of $\bar{B}$.

Let $t_1=t,\dots,t_n$ generate $\M$ minimally. We may replace $S$ by
the germ of curve defined by $t_2,\dots,t_n$. Thus $t$ is a
uniformizing parameter of $A$.

Since the map germ of $D\ar S$ has an ordinary quadratic singularity
at $P$, there are regular parameters $\bar{x}_1$,$\bar{x}_2$ of $D$
around $P$ such that $t\mapsto\bar{x}_1\bar{x}_2$.

So now we have reduced to the following. The completion of the local
ring of $S$ at the image of $P$ may be writen as $A[|t|]$ for some
power series ring $A$. The completion of the local ring of $X$ (resp.
$D$) at $P$ is of the form $B=A[|t,x_1,x_2|]$ (resp. $B/(t-x_1x_2))$. Hence
$S'$ is represented by the ideal $(t,x_1,x_2)\subset B$. The diagonal and
the fibre product $X\times_S{X}$ are represented by $(x_1-
x_1',x_2-x_2') \subset{A}[|t,x_1,x_2,x_1',x_2'|]$. The blowup
$X_2\ar{X}\times_S{X}$ is given by the inclusion
${A}[|t,x_1,x_2,x_1',x_2'|]\subset{A}[|t,x_1,x_2,x_1',u|]$ defined by
$x_2'=x_2+u(x_1'-x_1)$. Restriction over $S'$ therefore takes on the
form, $A\ar{A}[|x_1',x_2'|]\subset{A}[|x_1',u|]$, with $D'$ defined by
$u$. \qed

\bp
\label{reduced}
Let $Y$ be a smooth, projective surface and let $D$ be an ample
divisor on $Y$. Fix $n\!\in\!\{1,\dots,6\}$. Then there exists an
integer $r_0$ such that, for all $r\geq r_0$, for all linear
subsystems $S$ of $|rD|$ of dimension $n$ in an open dense subset
of the appropriate grassmannian, the following holds:
\bi \item[] $\Sigma((2^{[n]});S)$ is finite, reduced,
and for $(C,y_1,\dots,y_n)\!\in\! \Sigma((2^{[n]});S)$ we have that
$(y_1,\dots, y_n)$ is a singularity of one of the strict types
described in Prop.\ref{tipos}.
\ei\ep
\pf As observed in the proof of Prop.\ref{tipos}, ampleness ensures
that for any fixed sequence $\m=(m_1,\dots, m_n)$ of positive
integers there exists $r_0$ such that, for all $r\geq r_0$, and for
any sequence $(y_1,\dots,y_n)$ of distinct points in $Y$, the sheaf
$\M_{y_1}^{m_1}\cdots\M_{y_n}^{m_n} \otimes \O(rD)$ is generated by
global sections. It follows that distinct $y_i$'s impose independent
conditions to be a singularity of strict type $\m$ on the system $|rD|$
and in fact, $\Sigma(\m,rD)$ restricted to the complement of the
union of the exceptional divisors in $Y_n$ is a projective bundle. In
\cite{i} ((9.1),p. 417) it is shown the same is true over all of
$Y_n$ provided $\m$ satisfies the relaxed proximity inequalities
$m_{i}\geq- 1+m_{i+1}+\cdots+m_n$ for $i=1,\dots,n-1$.

As this is no longer the case for $\m=(2^{[n]}),\ n\geq3$, a direct
verification of smoothness is required.

At any rate, $\Sigma((2);S)$ and $\Sigma((2,2);S)$ are smooth for
all sufficiently ample complete system $S$ and remain smooth
upon replacing $S$ by a general subsystem by transversality of a
general translate \cite{kltransv}.

For $n\geq3$ we proceed by the following iteration argument.
Recall from Prop.\ref{step} that for any
$D\subset X\ar S$ as in \S\ref{setup}, we have
$$\Sigma((2^{[3]});D)=\Sigma((2);D'),$$
where $D'=(f_{3,2}^*(f_{2,2}^*D-2E_{1,2})-2E_{1,3})_{|\Sigma((2,2);D)}$.

If $S$ is a sufficiently ample complete system, one checks that $D'$
is regular. In fact, it is the total space of a family of
basepoint--free divisors in the fibres of $Y_3\ar{Y_2}$. Indeed, let
$Y'\ar{Y}$ be the blowup at $y_1\!\in\!{Y}$ and let $Y''\ar{Y'}$ be the
blowup at $y_2\!\in\!{Y'}$. Let $y_3\!\in\!{Y''}$. Let $\L$ be an ample line
bundle over $Y$. Then
$$H^1(Y'',\L^{\otimes{r}}\otimes{\O_{Y''}(-2E_{y_2})}\otimes {\O_{Y'}(
2E_{y_1})}\otimes{\M_{y_3}})=0$$ for $r>>0$ because the sequence
$(2,2,1)$ satisfies the relaxed proximity inequalities. Hence
Prop.\ref{iter} implies that $\Sigma((2^{[3]},D)$ is regular at any
$(C,y_1,$ $y_2,$ $y_3)$ such that $y_3$ is an {\em odp} of
$C-2E_{y_1}-2E_{y_2}$. Now, if $y_3$ were a triple point (allowed if
$n=6$), then we would certainly have $y_3$ not infinitely near ${y_2}$.
Let $\pi$ be the involution of $X_2 \times_{X} X_2$ (so that
$p_{3,2}\pi=p_{3,1}$). It lifts to an involution of $X_3$ that leaves
$\Sigma((2^{[3]},D)$ invariant. Since $\pi$ maps $(C,y_1,y_2,y_3)$ to
$(C,y_1,y_3,y_2)$, we get regularity at the latter as well. The same
argument yields regularity of $\Sigma((2^{[n]},D)$ for $n=4,5,6$ and
$S$ generic, $\infty^n$. For instance, to show $\Sigma((2^{[6]},D)$
is regular at $(C,y_1,\dots,y_6)$ such that $y_1$ is of strict type
$(3(2))$ and $y_2$ is the double point infinitely near (cf.\ref{tritwo}),
we argue by regularity of $\Sigma((2^{[2]},D)$ at $(y_1,y_2)$ and apply
iteration, observing that $(y_3,\dots,y_6)$ are all ordinary
quadratic singularities. If $y_2$ were the intersection of the
exceptional line and the smooth branch, then $y_3$ must be the double
point infinitely near to $y_1$. In this case apply first a
permutation and argue as before. \qed
\section{Applications}
Here are two situations we may apply the above constructions to.

\n1:\underline{Linear systems}. Let $Y$ be a smooth projective surface, let
$\M$ be an invertible \O$_Y-$module and let $V\subset H^0(Y,\M)$ be a subspace.
Set $S=\p{}(V^*), X=S\times Y$ and let $f:X\ar S$ be the projection. Then
$\L=\M\otimes \O(1)$ has a section defining the universal divisor $D$ of the
linear system parametrized by $S$. We also write in this case,
$\Sigma(\m;S):=\Sigma(\m;D)$.
\vskip10pt
\n2:\underline{Hypersurfaces}. Let $S=Gr(2,N)$ be the Grassmann variety of
planes in \p{N}, with tautological quotient sheaf
$\O^{N+1}\surj \Q$, where rank$\Q$=3. Let $X=\p{}(\Q)\subset
S\times\p N$ be the universal plane in \p N. Set
$\L=\O\!_{\pp{}({\cal Q})}(d)$ and let $D\subset X$ be defined by a form
of degree $d$. Thus the fibre of $D$ over $s\!\in\! S$ is the
intersection of a fixed hypersurface with the plane $s$
represents.

Using Prop. \ref{reduced} we get the following formulas
for the number  $tg_n$ of $n-$nodal curves in an $\infty^n$
family of curves, for $n\!\in\!\{1,\dots,6\}$.

\bp
\label{formulas}
Fix $n\!\in\!\{1,\dots,6\}$.
Let $D\subset X\ar S$ be a family of curves in a smooth family
of surfaces of dimension $n$. Assume $\Sigma((2^{[n]};D)$ is
reduced and receives contributions only from the strict types
described in Prop.\ref{tipos}. Then we have:
\bi
 \item[] $tg_n=(\#\Sigma((2^{[n]});S)/n!\hbox{\quad for }n\!\in\!\{1,2,3\};$
 \item[] $tg_4=(\#\Sigma((2^{[4]});S)-6\#\Sigma((3);S))/4!;$
 \item[] $tg_5:=(\#\Sigma((2^{[5]});S)-30\#\Sigma((3,2);S))/5!;$
 \item[] $tg_6:=(\#\Sigma((2^{[6]});S)-30\#\Sigma((3(2));S)-
 90 \#\Sigma((3,2,2);S))/6!.$
\ei\ep
\pf
Let us explain for instance the coefficient 90 appearing in
the formula for $tg_6$.
Pick $(C,z_1,z_2,z_3)$ in $\Sigma((3,2,2);S)$. Here $C$ is a curve in the
system $S$ and $(z_1,z_2,z_3)$ is a singularity of strict type $(3,2,2)$.
Let $z_{11},z_{12},z_{13}$ be the branches over $z_1$.
It gives rise to the following list of singularities $(y_1,\dots,y_6)$
of weak type $2^{[6]}$ on $C$:
\ben
\item[] $y_1=z_1$ and $(y_2,\dots,y_6)=$ any permutation of
$\{z_2,z_3,z_{11},z_{12},z_{13}\}$\vskip1pt\hfill SUBTOTAL:~120.

\item[] $y_1=z_i ,~y_2=z_1$ and $(y_3,\dots,y_6)=$ any permutation of
$z_j,z_{11},z_{12},z_{13}$ with $\{i,j\}=\{2,3\}$\hfill SUBTOTAL:~48.

\item[] $y_1\!=\!z_i,y_2\!=\!z_j,y_3\!=\!z_1$ and $(y_4,y_5,y_6)=$ any
permutation of $z_{11},z_{12},z_{13}$ with $\{i,j\}=\{2,3\}$
						\hfill SUBTOTAL:~12.
\en
The factor 180/2 comes from the fact that $(C,z_1,z_2,z_3)$ and
$(C,z_1,z_3,z_2)$ yield the same contributions to $\Sigma((2^{[6]});S)$
 \qed  \vskip10pt

Using the formula (\ref{dejonq}) in Prop.\ref{step} the $rhs$ can
be computed in terms of Chern classes for each of the two
situtions envisaged above. We've made extensive use of {\sc maple\cite{maple}\
\& schubert\cite{schub}}). See the appendix for the
computations.
\section{Surfaces}
For the case of linear systems on a surface $Y$, setting for short,
$$\ba{ll}
c_2=degree(c_2\Omega^1_Y),\ &k_1=degree(c_1\Omega^1_Yc_1\L),\\
k_2=degree((c_1\Omega^1_Y)^2),\ &d=degree((c_1\L)^2).
\ea$$
we get from (\ref{sh}), (\ref{formulas}) and (\ref{dejonq}),
\vskip8pt
\footnotesize
\n$tg_1:=3 d\!+\!2 k_1\!+\!c_2;
$\vskip8pt\n$
tg_2:=(tg_1 (\!-\!7\!+\!3 d\!+\!2 k_1\!+\!c_2)\!-\!6 k_2\!-\!25 k_1\!-
\!21 d)/2
;$\vskip8pt
\n$tg_3:=(2tg_2 (\!-\!14\!+\!3 d\!+\!2 k_1\!+\!c_2)\!+\!tg_1 (\!-
\!6 k_2\!-\!25 k_1\!-\!21 d\!+\!40)\!+\! (\!-\!6 k_2\!-\!25
k_1\!-\!21 d) c_2\!-\!63 d^2\!+\!(\!-\!18 k_2\!-\!117
k_1\!+\!672) d\!-\!50 k_1^2\!+\! (\!-\!12 k_2\!+\!950)
k_1\!+\!292 k_2)/6;$\vskip8pt
\n$tg_4= (81d^4\!+\!(216k_1 \!+\! 108c_2 \!-\! 2268)d^3\!+\!
(54c_2^2 \!+\! (216k_1 \!-\! 1890)c_2 \!-\! 324k_2 \!+\! 21852
\!-\! 5130k_1 \!+\! 216k_1^2)d^2 \!+\! (12c_2^3 \!+\! ( \!-\!
504 \!+\! 72k_1)c_2^2 \!+\! ( \!-\! 216k_2 \!+\! 8940 \!+\!
144k_1^2 \!-\! 2916k_1)c_2 \!-\! 3816k_1^2 \!+\! 39780k_1 \!+\!
96k_1^3 \!+\! 6024k_2 \!-\! 72360 \!-\! 432k_1k_2)d \!+\! c_2^4
\!+\! ( \!-\! 42 \!+\! 8k_1)c_2^3 \!+\! ( \!-\! 402k_1 \!-\!
36k_2 \!+\! 24k_1^2 \!+\! 699)c_2^2 \!+\! ( \!-\! 3888 \!-\!
144k_1k_2 \!+\! 1756k_2 \!+\! 9046k_1 \!-\! 1104k_1^2 \!+\!
32k_1^3)c_2 \!-\! 144k_1^2k_2 \!+\! 16k_1^4 \!+\! 108k_2^2 \!+\!
4412k_1k_2 \!-\! 936k_1^3 \!+\! 17171k_1^2 \!-\! 28842k_2 \!-\!
95670k_1)/24$;
\vskip8pt\n$
tg_5 = 81/40 d^5\!+\!(27/8 c_2\!+\!27/4 k_1\!-\!189/2)
d^4\!+\!(9/4 c_2^2\!+\!(\!-\!441/4\!+\!9 k_1)c_2\!+\! 9 k_1^2\!-
\!27/2 k_2\!-\!1107/4 k_1\!+\!3393/2) d^3\!+\!(3/4 c_2^3\!+\!(\!-
\!189/4\!+\!9/2 k_1) c_2^2\!+\! (9 k_1^2\!-\!981/4
k_1\!+\!2469/2\!-\!27/2 k_2) c_2\!-\!27 k_1 k_2\!+\!6 k_1^3\!-
\!603/2 k_1^2\!-\! 13875\!+\!471k_2\!+\!8463/2 k_1) d^2\!+\!(1/8
c_2^4\!+\!(\!-\!35/4\!+\! k_1) c_2^3\!+\!(3 k_1^2\!-\! 285/4
k_1\!+\!2207/8\!-\!9/2 k_2) c_2^2\!+\!(4 k_1^3\!-\!4789\!-\!18
k_1 k_2\!-\!180 k_1^2\!+\! 565/2 k_2\!+\!8589/4 k_1) c_2\!-\!145
k_1^3\!-\!22445/4 k_2\!+\!27403/8 k_1^2\!+\!2 k_1 ^4\!+\!27/2
k_2^2\!+\!1355/2 k_1 k_2\!-\!111959/4 k_1\!+\!217728/5\!-\!18
k_1^2 k_2) d\!+\! 1/120 c_2^5\!+\!(1/12 k_1\!-\!7/12)
c_2^4\!+\!(141/8\!+\!1/3 k_1^2\!-\!1/2 k_2\!-\!27/4 k_1)
c_2^3\!+\!(251/6 k_2\!-\!53/2 k_1^2\!-\!3 k_1 k_2\!+\!2/3
k_1^3\!-\!485/2\!+\!1547/6 k_1) c_2^2 \!+\!(\!-\!17881/12
k_2\!+\!3516/5\!+\!1229/6 k_1 k_2\!-\!68137/12 k_1\!-\!131/3
k_1^3\!+\!9/2 k_2^2\!+\!21551/24 k_1^2\!+\!2/3 k_1^4\!-\!6 k_1^2
k_2) c_2\!+\!727/3 k_1^2 k_2\!-\! 188k_2^2\!-\!8827/2 k_1
k_2\!+\!321882/5 k_1\!+\!9 k_2^2 k_1\!+\!22695 k_2\!+\! 10867/12
k_1^3\!-\! 26189/2 k_1^2\!-\!4 k_1^3 k_2\!+\!4/15 k_1^5\!-\!26
k_1^4;$
\vskip8pt\n$
tg_6=
( 81/80 ) d^6\!+\!( 81/40 c_2\!-\!567/8\!+\!81/20 k_1 ) d^5\!+\!(
27/16 c_2^2\!+\!(27/4 k_1\!-\!1701/16) c_2\!-\!81/8
k_2\!+\!8109/4\!+\!27/4 k_1^2\!-\!4077/16 k_1 ) d^4\!+\! ( 3/4
c_2^3\!+\!(9/2 k_1\!-\!63) c_2^2\!+\!(8523/4\!-\!27/2 k_2\!+\!9
k_1^2\!-\!1233/4 k_1) c_2\!+\! 1131/2 k_2\!+\!6 k_1^3\!-
\!29601\!-\!27 k_1 k_2\!-\!729/2 k_1^2\!+\!25671/4 k_1 ) d^3\!+\!
( 3/16 c_2^4\!+\!(3/2 k_1\!-\!147/8) c_2^3\!+\!(12909/16\!-\!27/4
k_2\!+\!9/2 k_1^2\!-\!1107/8 k_1) c_2^2\!+\!(2073/4 k_2\!-
\!76959/4\!-\!27 k_1 k_2\!+\!41493/8 k_1\!+\!6 k_1^3\!-\!333
k_1^2) c_2\!+\! 3 k_1^4\!+\!81/4 k_2^2\!-\!27 k_1^2 k_2\!-
\!96699/8 k_2\!-\!519/2 k_1^3\!+\!1102009/5\!+\! 119961/16
k_1^2\!-\!639927/8 k_1\!+\!4821/4 k_1 k_2 ) d^2\!+\! ( 1/40
c_2^5\!+\!(1/4 k_1\!-\!21/8) c_2^4\!+\!(\!-\!3/2
k_2\!+\!3071/24\!-\!109/4 k_1\!+\!k_1^2) c_2^3 \!+\!(\!-\!201/2
k_1^2\!+\!157 k_2\!+\!2 k_1^3\!-\!29213/8\!-\!9 k_1
k_2\!+\!5421/4 k_1) c_2^2\!+\! (\!-\!26787/4 k_2\!+\!648997/10\!-
\!74149/2 k_1\!-\!159 k_1^3\!+\!1481/2 k_1 k_2\!+\!27/2
k_2^2\!+\! 32959/8 k_1^2\!+\!2 k_1^4\!-\!18 k_1^2 k_2)
c_2\!+\!853 k_1^2 k_2\!-\!18481 k_1 k_2\!-\! 1317/2 k_2^2\!+\!27
k_2^2 k_1\!+\!1401361/12 k_2\!+\!28988249/60 k_1\!+\!46109/12
k_1^3\!-\! 12 k_1^3 k_2\!+\!4/5 k_1^5\!-\!668388\!-\!554465/8
k_1^2\!-\!92 k_1^4 ) d\!+\!1/720 c_2^6\!+\!(\!-\!7/48\!+\! 1/60
k_1) c_2^5\!+\!(\!-\!1/8 k_2\!+\!1/12 k_1^2\!-\!95/48
k_1\!+\!331/48) c_2^4 \!+\!(\!-\!k_1 k_2\!-\!10 k_1^2\!+\!8147/72
k_1\!-\!8095/48\!+\!565/36 k_2\!+\!2/9 k_1^3) c_2^3\!+\! (\!-
\!145/6 k_1^3\!+\!15347/10\!+\!1355/12 k_1 k_2\!-\!3 k_1^2
k_2\!+\!1/3 k_1^4\!+\!9/4 k_2^2\!-\! 190339/48 k_1\!+\!26519/48
k_1^2\!-\!10891/12 k_2) c_2^2\!+\!(\!-\!4 k_1^3 k_2\!-\!85/3
k_1^4\!+\! 4291/4 k_1^3\!+\!9 k_2^2 k_1\!+\!10998\!-\!815/4
k_2^2\!-\!807341/48 k_1^2\!+\!790/3 k_1^2 k_2\!+\! 4/15 k_1^5\!-
\!62339/12 k_1 k_2\!+\!691883/24 k_2\!+\!$ $10672201/120 k_1)
c_2\!-\! 311237/16 k_1^3\!-\!9/2 k_2^3\!+\!4/45
k_1^6\!+\!7001519/72 k_1 k_2\!-\!2 k_1^4 k_2\!-\! 1855/4
k_2^2 k_1\!+\!9 k_1^2 k_2^2\!+\!1805/9 k_1^3 k_2\!-
\!1080646 k_1\!+\! 86753363/360 k_1^2\!+\!200477/36
k_2^2\!+\!26297/36 k_1^4\!-\!13 k_1^5\!-\! 55951/8 k_1^2 k_2\!-
\!2567321/6 k_2.$
\normalsize
\bex{\bf  $Y=\p{2}$.}
\label{p2}
We make the substitutions,\ex
\centerline{$c_2=3,d=m^2,k_1=-3 m,k_2=9.$}
 \bs{ $n=4.$}
 The expression for $tg_4$ above reduces to\es
\footnotesize $$tg_{4,\p{2}}(m)=\!-\!
8865\!+\!18057/4 m\!+\!37881/8 m^2\!-\!2529 m^3\!-\!642 m^4\!+\!1809/4 m^5
\!-\!27m^7\!+\!27/8 m^8.$$
\n\normalsize Setting $m=4$ we get
$666=126+540$
for the number of 4--nodal quartics through 10 general points.
Indeed,  a plane quartic with 4 nodes splits as a
union of 2 conics, 126 of which pass through 10 points, or of a
singular cubic and a line through $10$ points.

\bs{$n=5.$}
We find,\es
\n\footnotesize
$tg_{5,\p{2}}(m)=
81/40 m^{10}\!-\!81/4 m^9\!-\!27/8 m^8\!+\!
2349/4 m^7\!-\!1044 m^6\!-\!127071/20 m^5\!+\!128859/8
m^4\!+\!59097/2 m^3-3528381/40 m^2-946929/20 m\!+\!153513$.
\normalsize
\vskip8pt\n
Setting $m=4$ and picking a system of quartics through 9
general points, we do get the right answer,
378=\bin{9}{5}$\times$3.
Indeed, a plane quartic with 5 nodes can only be a union of a
conic and line pair: hence \bin{9}{5}
for the choice of 5 points determining a conic, times the
number 3 of line pairs through the 4 remaining points...

\bs{ $n=6.$} We have,\par
\n\footnotesize
$tg_{6,\p{2}}(m)=$
$81/80 m^{12}$ - $243/20 m^{11}$ $-$ $81/20 m^{10}$ + $8667/16
m^9$ $-$ $9297/8 m^8 $ $-$ $47727/5 m^7$ + $2458629/ 80 m^6$ +
$3243249/40 m^5$ $-$ $6577679/20 m^4$ $-$ $ 25387481/80 m^3$ +
$6352577/4 m^2$ + $8290623/20 m $ $-$ $2699706.$
\normalsize
\vskip8pt\n
Again setting $m=4$, we find 105 for the number of 6--nodal
quartics through 8
general points: the configurations must consist of 4 lines.
\bss
Setting $m=5$, we can find the number of {\em irreducible} rational
plane quintic curves through 14 general points. This is
$tg_{6,\p{2}}(5)-\bin{14}{2}tg_{2,\p{2}}(4)-\bin{14}{5}$
$= 109781-20475-2002=87304$.
The corrections are due to the reducible 6--nodal quintics:
either line+binodal quartic or conic+cubic.
\ess\es
\bex{\bf$Y=\q$.} For a system of curves of type $(m_1,m_2)$, we
set\newline
\centerline{$c_2=4,k_2=8,k_1=-2(m_1+m_2),d=2 m_1m_2$.}\ex
\bs{$n=4.$}
We get,
\vskip8pt\n\footnotesize
$tg_{4,\qp}(m_1,m_2) =
(32/3\!-\!64 m_2\!+\!144 m_2^2\!-\!144 m_2^3\!+\!54 m_2^4)
m_1^4\!+\!( 808/3\!-\!3112/3 m_2\!+\!1230 m_2^2\!-\!324 m_2^3\!-
\! 144 m_2^4) m_1^3\!+\! (11987/6\!-\!3494 m_2\!-\!2 m_2^2\!+\!
1230 m_2^3 \!+\! 144 m_2^4 ) m_1^2 \!+\!(17359/6\!+\!11333/3
m_2\!-\!3494 m_2^2\!-\!3112/3 m_2^3 \!-\!64 m_2^4) m_1 \!-
\!7460\!+\!17359/6 m_2\!+\!11987/6 m_2^2\!+\!808/3 m_2^3  \!+\!
32/3 m_2^4.
$\vskip10pt
\normalsize
\bss
If $m_1=m_2=2$, it checks with the number $6$ of configurations
of $4$ lines in the system $(2,2)$ through  4 general points on a
quadric. Indeed, since $p_a=1$, the curve splits in one of the types:
$(1,1)+(1,1)$, $(2,0)+(0,2)$, $(2,1)+(0,1)$ or
$(1,2)+(1,0)$. The latter two cases consist of the union of a
twisted cubic and a bi-secant line, hence get for free two nodes
due to the intersections. In order to present $4$ nodes, the
twisted cubic must split further. One easily sees that the only
possibility is indeed a configuration $(2,0)+(0,2)$ of $4$ lines.
 We may assume no $2$ of the
$4$ points are on a ruling. Label the points $1,2$ so that the
lines composing the curve $(2,2)$ through them are both of system
$(1,0)$; this forces the other $2$ lines to be of the opposite
system $(0,1)$. Thus, the choice of $1,2$ completely determines
the solution, hence \bin{4}{2}.
\ess
\bss
For $(m_1,m_2)=$ $(2,3)$, we find
$tg_{4,\mbox{\footnotesize\q}}(2,3) =133$. As $p_a=2$,
we obtain again reducible configurations.
Notice the system $|(2,3)|$ is $\infty^{11}$.
Let the $\infty^4$ subsystem be defined by imposing 7 points.
Possible splitting types?
(i)$(2,0)+(0,3)$ is $\infty^{5}$, too small.
(ii)$(2,1)+(0,2)$ is $\infty^{7}$; 4 nodes due to intersection,
 $\bin{7}{2}=21$ choices for configuration consisting of twisted
 cubic$\in\!|(2,1)|$ through $5$ points and line pairs$\in\!| (0,2)|$
 through $2$ points.
 SUBTOTAL:~ 21.
(iii)$(1,1)+(1,2)$ is $\infty^{3+5}$; $3$ nodes due to
intersection, hence need additional node for either $(1,1)$ or
$(1,2)$ component. If the new node is on $(1,1)$, this curve must
be a line pair; make it pass through 2 of the points  ($\bin{7}{2}$
choices for these) $\times$ 2 (=number of such line pairs for
each choice of 2 points), total  42. One takes the $(1,2)-$component
through the remaining 5 points, unique choice. SUBTOTAL:~ 42.
If the new node is on a $(1,2)-$curve, this must split as $(0,1)
+(1,1)$, so the  actual solutions are of
the form $(1,1)+(0,1)+(1,1)$; if the  $7\up{th}$ point is on the line,
the remaining 6 will be on $\bin{6}{3}/2$ conic pairs. SUBTOTAL:~  70.
(iv)$(2,2)+(0,1)$ has  2  nodes due to intersection, hence need two
additional nodes for $(2,2)-$component; now
if a $(2,2)-$curve acquires  2  double points, it splits as
$(2,1)+(0,1)$ or $(1,2)+(1,0)$ or $(1,1)+(1,1)$; these have already
been accounted for! Thus it all happily adds up to the
right TOTAL:~ 133.
\ess
\bss
For $(m_1,m_2)=$ $(2,4)$, we find
$tg_{4,\mbox{\footnotesize\q}}(2,4) = $ 1261.
The system $|(2,4)|$ is $\infty^{14}$. We impose 10 points
to select an $\infty^{4}$ subsystem.
Possible splitting types?
(i) $(2,3)+(0,1)$ is $\infty^{12}$; $2$ nodes due to intersection.
Impose $2$ new nodes for $(2,3)-$component; there are
$tg_{2,\mbox{\footnotesize\q}}(2, 3)=105$ through each of the 10
choices of  9 points. Notice that
among these  1050 curves
there are  90 in $|(2,2)+(0,2)|$. These will
be accounted for separately  below.
SUBTOTAL:~  960.
(ii) $(2,2)+(0,2)$ is $\infty^{10}$;  4 nodes due to
intersection; $\bin{10}{2}=45$ choices for 2 points determining
a line pair in the system $|(0,2)|$, the remaining 8 points
singling out a member in $|(2,2)|$.
SUBTOTAL:~ 45.
(iii) $(2,1)+(0,3)$ is $\infty^{8}$: too small! Similarly for
$(2,0)+(0,4)$.
(iv) $(1,4)+(1,0)$ is $\infty^{10}$;  4 nodes due to
intersection;  10 choices for the point determining the
component $(1,0)$.
SUBTOTAL:~  10.
(v) $(1,3)+(1,1)$ is $\infty^{10}$;  4 nodes due to
intersection; $\bin{10}{3}$ choices for  3 points determining
a conic while the 7 other points determine the
component $(1,3)$.
SUBTOTAL:~  120.
(vi) $(1,2)+(1,2)$ is $\infty^{10}$;  4 nodes due to
intersection; $\bin{10}{5}$ choices for  5 points determining a
twisted cubic $(1,2)$
SUBTOTAL:~  126.
It gives the expected TOTAL:~ 1261.
\ess
\bss{\bf Irreducible rational curves with $p_a=4$ on \q.}\label{d33}
We may compute the number {\bf3510} of
irreducible rational curves of type $|(3,3)|$ passing
through $11=15-4$ general points. We subtract from
$tg_{4,\mbox{\footnotesize\q}}(3,3) =4115$,
the contributions given by:
(i) \big(nodal $(3,2)$ through 10 points + $(0,1)$ through
the 11\up{th}\big):
	$20\times11=220$;
(ii) \big(nodal $(2,3)$ through 10 points + $(1,0)$\big): 220;
(iii) \big($(2,2)$ through 8 points + $(1,1)$ through 3
others\big): \bin{11}{8}=165.
(Note that $(3,1)+(0,2)$ is $\infty^{9}-$ too small.)
\ess\bss\label{d25}
Reasoning as above, we also find the number {\bf3684} of
irreducible rational curves in the system $(2,5)$ passing
through $13=17-4$ general points. This is
\big($tg_{4,\mbox{\footnotesize\q}}(2,5) =7038$\big) minus
\big($(tg_{2,\mbox{\footnotesize\q}}(2,4) =252)\times13$
due to binodal $(2,4)+(0,1)$\big)
minus \big(\bin{13}{11}=78 due to curves $(2,3)$
through 11 points + $(0,2)$ through 2 others\big).
\ess\es
\bs{$n=5.$} The first interesting check is provided by
the system $|(3,3)|$ on $\q$. We
find $tg_{5,\mbox{\footnotesize\q}}(3,3) =3702$. Here we have
$p_a=4$, hence imposing 5 nodes will force again reducible curves.
Fix 10 points in general position to
define an $\infty^5$ subsystem of $|(3,3)|$.
Possible splitting types?
(i)$(3,1)+2(0,1)$ is $\infty^{7+1} $ and $(3,1)+(0,2)$ is
$\infty^{7+2}$, both too small. (ii)$(3,2)+(0,1)$:
$\infty^{12}$; there are 3 nodes due to intersection. Look at
members of$|(3,2)|$ through 9 points and with 2 additional nodes: we find
$tg_{2,\mbox{\footnotesize\q}}(3,2)=105$. Among these, 9 split further
as $(2,2)+(1,0)$ and will be accounted for separately in (iv).
Since there are 10 choices for the  9 points,
we have the SUBTOTAL:~ 960. (iii)$(2,3)+(1,0)$: just as in (ii),
SUBTOTAL:~960. (iv)$(2,2)+(1,0)+(0,1)$: there are
$\bin{10}{2}=45$ times 2 for choices of points and system of
line through them. SUBTOTAL:~90.
(v)$(2,2)+(1,1)$: we have 4 nodes due to intersection. When the
aditional node is on the $(2,2)$ component which passes through
7 points, we find
$tg_{1,\mbox{\footnotesize\q}}(2,2)=12$, times $\bin{10}{3}$ obtaining
the SUBTOTAL:~1440.
If the additional node be on $(1,1)$, the type becomes
$(2,2)+(1,0)+(0,1)$, already accounted for in (iv) above.
(vi)$(2,1)+(1,2)$:there are 5 nodes due to intersection; contributes
$\bin{10}{5}$,
SUBTOTAL:~ 252,
fortunately totaling  3702.
\bss{How about the irreducible rational curves with $p_a=5$ on \q?}
The possible bidegrees are $(2,6),(6,2)$. One expects finitely many of
these passing through 15 points. However we notice that {\em any}
subsystem $S\subset|(2,6)|$ of codimension 15 meets the
family of curves of type $(2,4)+2(0,1)$. Since
these present a nonreduced component, therefore
$\Sigma((2^{[5]});S)$ contains components of wrong dimension
(cf.\ref{nonred}), so that the formula is {\em not} applicable to the
present case. It would be nice to compute the equivalence of
these bad components.

\ess\es
\bs{$n=6.$}
\bss We look again at the system $|(3,3)|$ on $\q$. We
find  $tg_{6,\qp}$ $(3,3)=$ 2224.
Fix 9 points in general position to
define an $\infty^6$ subsystem.
Possible splitting types?
(i)$(3,1)+2(0,1):\ \infty^{7+1}$, too small.
(ii)\big(nodal$(2,2)$\big) $+(1,0)+(0,1)$:
$ 12\times\bin{9}{7}\times2$. SUBTOTAL:~864.
(iii)$(3,1)+(0,2):\infty^{7+2}$; contributes $\bin{9}{7}.$ SUBTOTAL:~ 36.
(iv)$(1,3)+(2,0)$: SUBTOTAL:~ 36.
(iv)
$(2,1)+(1,1)+(0,1):\infty^{5+3+1}$; contributes
$\bin{9}{5}\times\bin{4}{3}.$ SUBTOTAL:~        504.
(v)$
(1,2)+(1,1)+(1,0)$. SUBTOTAL:~ 504.
For several days, we had found only these 1944. The 280 then missing
were pointed out to me (after a lunch break at the MSRI) by
Enrique Arrondo: $\bin{9}{3}\times\bin{6}{3}/6=280$ curves
of the form $(1,1)+(1,1)+(1,1)$!!!
\ess
\bss{\bf Irreducible rational curves of
bidegree $(3,4)$}\label{d34} passing through $13=19-6$ general
points: {\bf90508}. We subtract from
$tg_{6,\mbox{\footnotesize\q}}(3,4) =122865 $, the contributions
given by: (i)\big(trinodal $(3,3)$ through 12 points\big) +
\big($(0,1)$ through the 13\up{th}\big): $1944\times13=25272$;
(ii)\big(nodal $(2,3)$ through 10
points\big) + \big($(1,1) $ through 3 others\big):
20$\times$\bin{13}{3}=5720; (iii) \big($(2,2)$ through 8 points\big) +
\big($(1,2)$ through 5 others\big): \bin{13}{8}= 1287; (iv)\big($(3,2)$
through 11 points\big) + \big($(0,2)$ through 2 others\big):
\bin{13}{2}=78.
\ess\es

\bex{\bf Del-Pezzo surface:}
 $Y=\p{2}$ blown up at  5
points, imbedded in \p{4} as a $(2,2)$ intersection by the system
of plane cubics through the 5 points. There are  \bf40 \rm
fourfold tangent hyperplanes. Indeed, label the points $\{1,\dots,5\}$;
draw the lines $\overline{12},\ \overline{15},\ \overline{34}$;
let $a=\overline{12} \cap \overline{34},\ b= \overline{15} \cap
\overline{34}.$ Note $1$ is double on $\overline{12} +
\overline{15} + \overline{34}$. After blowing up, the hyperplane
system $|3L-e_1-\cdots-e_5|$ will contain the curve $e_1 +
\overline{12}' + \overline{15}' + \overline{34}'$ (the $'$
denoting strict transform).
It presents the 2 double points $a',b'$ and two others on $e_1$.
The number of such configurations
can be counted as  5 choices for the point labeled 1,
times \bin{4}{2} choices for
$\overline{12}, \overline{15}$, totaling 30. In addition to these
configurations of lines, we may also take the conic $c$ and a
line through a pair of the points, say $\overline{12}$; then we
get the hyperplane section $c' + \overline{12}' +e_1 +e_2$. This
gives  10 more, totaling 40, as predicted by the formula.
\ex
\bex{\bf Surfaces of degree $9$ in \p{4}.}
 Substituting
$$[d=9,k_1=2p_a-11,k_2=6\chi-5p_a+23, c_2=12\chi-k_2]$$ in $tg_4$
with the list of possible pairs (cf. \cite{aure})
$[p_a=\hbox{sectional genus};\chi=(c_2+k_2)/12]$ yields the
table,
{\footnotesize
$$\ba{cccccccc}
[6;1]& [7;1]& [7;2]& [8;2]& [8;3]& [9;4]&[10;5]& [12;9]\\
\bf15645 &\bf 57162 &\bf 107646 &\bf 248671 &\bf 388846 & \bf
1022595 &\bf 2222868 &\bf 10957224\rm.
\ea$$}\ex
\bex{\bf$K3-$surfaces.}\label{k3} Let $Y$ be embedded by a complete
 system $|C|$ of curves of genus
$n\!\in\!\{3,4,5,6\}$. We have $2n-2=C\cdot(C+K_Y)$.
Substituting
$[d=2n-2,k_1=0,k_2=0, c_2=24]$ in $tg_n$
we find
$$\ba{rcccc}
n:  &3&4&5&6\\
tg_n:  &\bf3200&\bf 25650 &\bf 176256 &\bf 1073720\rm
\ea$$
For $n\!\in\!\{4,5,6\}$, the values given above for $tg_n$ are
smaller then those predicted by a formula Ciliberto and Lopez
(priv. communication) obtained by a degeneration argument.

A related development is the work of Manoil \cite{m}, where he addresses the
question of existence of rational points on $K3-$surfaces defined over
a number field. He proves the existence of curves of geometric genus
$\leq1$ for a certain class of surfaces by counting singular curves.
\ex
\bex{\bf Abelian surfaces }
$Y\subset \p{4}$. Here we find the number \bf150 \rm
of $4-$fold tangent hyperplanes. It might be more than just a
coincidence the fact that the contribution from $\#\Sigma(3;S)$
is also $=150$, suspiciously a factor of the number $15,000$ of
symmetries of the Horrocks-Mumford bundle, a generic section of
which is known to vanish precisely on $Y\dots$

The following comments were kindly communicated by Chad Schoen.

Let $Y$ be an Abelian surface with a polarization of type
$(1,5)$. Any Horrocks-Mumford Abelian surface is of this type.
The converse is almost true. I believe that any simple Abelian
surface with a $(1,5)$ polarization is a Horrocks-Mumford Abelian
surface. Let $N$ be an invertible sheaf giving the $(1,5)$
polarization. A curve in $|N|$ has self-intersection $10$. This
is the degree of the normal sheaf which is also the dualizing
sheaf. Thus the arithmetic genus is $6$. If the curve is
irreducible and has $4$ nodes it's normalization has genus  2. If
$Y$ is ``general'' its Picard number is $1$ and any hyperplane
section must be irreducible. Let $C$ be such a $4-$nodal curve
and $\tilde{C}$ its normalization. There is an isogeny
$Jac(\tilde{C}) \ar Y$ taking $\tilde{C}$ to $C$. Again if $Y$
has Picard number $1$, there is no choice but for this map to
have degree  5. Now the degree  5 unramified covers of $Y$ are
parametrized by the subgroups of order   5 in the fundamental
group of $Y$. Write $L$ for this lattice and $L'\subset L$ for
the index  5 subgroup. Assuming that $Y$ has Picard number $1$,
the  5 fold cover $f:J\ar Y$ will be the Jacobian of a genus  2
curve  if and only if $J$ is principally polarized. This will
occur if and only if the pull back of the $(1,5)$ polarization on
$Y$ is  5 times a polarization on $J$. In terms of lattices and
the Riemann form associated to the polarization we have:
$$A:(1/5L)/L \times L/5L \ar(1/5\z)/\z=\z/5.$$ This alternating
form on the  5 torsion of $Y$ has a two dimensional radical--call
it $K$. ($K$=vectors in $(1/5)L/L$ which are orthogonal to the
whole space). Now the pull back to $J$ is divisible by  5 if and
only if the restriction of $A$ to  $(1/5)L'/L \times L'/5L \ar
(1/5)\z/\z$ is identically zero. This occurs exactly when $K$
lies in $(1/5)L'/L$. We can count all such $L'$. They are
hyperplanes in \p{3} containing a fixed \p{1} all over the field
$\z/5$. Thus the $L'$ 's are parametrized by \p{1}$(\z/5)$. There
are  6 possible $L'$ 's. Thus  6 possible $J$'s.  Finally we note
that translation by elements of $K=\z/5\times \z/5$ give
automorphisms of $Y$ preserving the $(1,5)$ polarization. This
gives  $6 \times25 =150 $ four-nodal hyperplane sections. There
are only  6 different  isomorphism classes of genus  2 curve
which occur as normalizations.

{\bf Question:} Inversion in the Abelian variety should also
preserve the polarization (I {\it (C. Schoen)} think).
How does this permute the 4--nodal hyperplane sections?
\ex
\section{Threefolds}
The same method yields the formula,
\vskip10pt\footnotesize
\n$
tg_{6,m}=(
m^{18}-12m^{17}+24m^{16}+155m^{15}-405m^{14}
+1082m^{13}- 18469m^{12}+66446m^{11}
 - 192307m^{10}+1242535m^{9}
-4049006m^{8}+11129818m^{7}- 53664614m^{6}+166756120m^{5}
 -415820104m^{4}+ 1293514896m^{3}- 2517392160m^{2}+1781049600m)/6!
$
\vskip10pt\n\normalsize
for the number of planes in \p{4} that are  6-fold tangent to a
hypersurface of degree $m$.

\subsection{Quartics.}
For $m=4$, the formula above gives \bf5600\rm. This can be verified by
the following direct calculation via the Fano variety $F$ (cf.\cite{ak}) of
$\infty^1$ lines contained in a 4\up{ic} threefold {\bf T}. Presently
the counting refers to the set
$$\{(\ell_1,\dots,\ell_4)\in
F^{\times4}|\exists\hbox{ plane }\pi\ s.t.\ \ell_1+\cdots+\ell_4=
\pi\cap {\bf T}\}$$
 of 4--tuples of coplanar
lines in that family.

Let $\cs_i\inj\O^{\oplus5}\surj\Q_i$ (rank $\Q_i= i+1$) denote the
tautological sequence over the Grassmann variety
$G_i:=Gr(i,4)$ of $i-$dimensional subspaces of \p{4}.
Go to the incidence variety
$\I:=\{(\ell,\pi)\!\in\! G_1\times{}G_2|\ell\subset\pi\}.$
It carries the diagram of locally free sheaves, (omitting
pullbacks)
\be\label{taut}\ba{ccccc}
\cs_2   & \inj & \cs_1      & \surj        &\M\\
||              &      &\injdown&    &\injdown\\
\cs_2   & \inj &\O^{\oplus5}&\surj &\Q_2\\
  &      & \surjdown&    &\surjdown\\
  &      &   \Q_1     & = & \Q_1
\ea\ee
The universal plane \p{}$(\Q_2)$ contains the total space $D$
of the family of intersections with the fixed 4\up{ic}
hypersurface. Our goal is to compute the intersection class
supported by
$$\I_3:=\{(\ell_1,\ell_2,\ell_3,\pi)\!\in\!{}\I\times_{G_2}\I
\times_{G_2}\I | D_\pi\geq\ell_1+\ell_2+\ell_3\}.$$
Set
$\I_1=\{(\ell,\pi) | \ell\subset \pi\cap D_\pi\}.$
This is
expressible as zeros of a section of a suitable bundle. Indeed,
up on $\p{}(\Q_2)_{|\I}$, we have the Cartier divisors $D_{|\I}$
and $\l_1:=\p{}(\Q_1)_{|\I}$. One checks that $\I_1$ is exactly the
locus in $\I$ where ``$\l_1 \subset D$'' holds along fibers.
Studying the natural diagram of  $\O_{\pp{}(\Q_2)_{|\I}}-$modules,

$$\ba{l}
\O\\ \downarrow\hskip.7cm\sear{s}\\
\O(D)\ar{}\O_{\l_1}(D)
\ea$$
one sees that the slant arrow $s$ vanishes on the fiber over
$(\ell,\pi)\!\in\! \I$ iff $\ell\subset\pi\cap D_{\pi}$. Let
$p:{\p{}(\Q_2)_{|\I}}\ar {\I}$ denote the structure map; it follows that
$\I_1$ is the scheme of zeros of the section $p_*s$ of the direct
image $sym_4Q_1$ of $\O_{\l_1}(D)=\O_{\l_1}(4)$. We obtain
$[\I_1]=c_5sym_4Q_1$. Pulling back  $D$ to $\I_1$ (and abusing
notation), it splits as $D=D_1+\l_1$, thus defining $D_1$.
Moreover, since $\p{}(\Q_1)$ is the divisor of zeros of a section
of $\O_{Q_2}(1)\otimes{}\M^*$, we have
$\O(D_1)=\O_{Q_2}(4)\otimes\O_{Q_2}(-1)\otimes{}\M$. We may ask
when does $D_1$ split further.
Go to $\I_1\times_{G_2} \I$. Set $\l_2=\I_1 \times_{G_2} \l_1$ and
define $\I_2$ by imposing the fibers of $D_1$ to contain a
2\up{nd} line.

As before, $\I_2$ is given by the vanishing of a section of the
pushforward of $\O_{\l_2}(D_1)$. Denoting by $_{(i)}$ the
pullback to $\I\times_{G_2}\I\cdots$ via $i$\up{th} projection, we
find $[\I_2]=c_4(\M_{(2)}\otimes{}sym_3\Q_{1(2)}).$ Similarly,
pulling back $D_1$ over $\I_2$ yields $D_1=D_2+\l_2$
and we get $[\I_3]=c_3(\M_{(3)}\otimes{}sym_2Q_{1(3)})$. See in the
Appendix a script for the actual computation using
\schub\cite{schub}. Observing that a 6--fold tangent plane $\pi$
to a 4\up{ic} hypersurface cuts 4 lines, the computation gives
134400/24=\bf5600 \rm as asserted.

\subsection{Quintics}
Recall that a general 5\up{ic}threefold {\bf T}$\subset\!\p{4}$
 contains 2,875 lines and 609,250 conics (cf. \cite{har}, \cite{katz}).

The plane through a conic counts as a  6--fold tangent since its
intersection wih {\bf T} splits as a $conic+cubic$,
thereby presenting 6 nodes.

Through each line, there are $\infty^2$ planes in \p{4}. The intersection
of any such plane with {\bf T} splits as $line+quartic$ thereby
counting as a  4--fold tangent.
The plane is a 6--fold tangent iff the residual plane quartic is
binodal.

Fix a line $\ell\subset{\bf T}$.
Let us find, among these $\infty^2$ residual plane quartic
curves the number of those with  2 double points. This requires
the computation of $\Sigma((2,2);D)$ for the family
$D\!\subset{}X\ar{}S$ of residual plane quartic we now describe.
Notation as in the previous example, let
$\cs_2\inj\O^{\oplus5}\surj\Q_2$ (rank $\Q_2= 3$) denote the
tautological sequence over the Grassmann variety
$G_2$ of planes in \p {4}. Let $G_{2,\ell}$ be the
Schubert subvariety of all 2-planes through a fixed line $\ell$.
Let $X=\p{}(\Q_2)_{|G_{2,\ell}}{\subset}
G_{2,\ell}\times\p 4$ be the restriction over $G_{2,\ell}$
of the universal plane in \p 4.
Restricting the sequence over $G_{2,\ell} $ yields an exact sequence,
(cf. top sequence in (\ref{taut}))
$\cs_2\inj\O^{\oplus3}=\cs_{1|\ell}\surj \M$, where $\M$ is a line
subbundle of $\Q_2$ with Chern class
$x:=$ $c_1\M=-c_1\cs=c_1\Q_2$.
Over $X$, we have the natural commutative diagram of maps of locally
free sheaves,
$$\ba{ccc}
\M&&\\
\injdown&\sear{}&\\
\Q_2&\surj{}&\O\!_{{\cal Q}_2}(1)
\ea$$
where the bottom line is the tautological 1--quotient on the
projective bundle\break {\bf Proj}$(Sym(\Q_2))$.
 One checks that $\ell':=G_{2,\ell}\times\ell$ is the divisor in
$\p{}(\Q_2)_\ell$
of zeros of the slant arrow $\M\ar\O(1)$. Therefore, setting
$y=c_1\O(1)$ we have $\O(\ell')=\O(y-x)$.

Now let $ D_{\mbox{\footnotesize\bf T}}
\subset\p{}({\Q_2})$ be the divisor defined by
intersection with {\bf T}, so that $\O(D_{\mbox{\footnotesize\bf T}}
=\O(5\cdot y)$.
Restriction over $G_{2,\ell}$ splits $D_{\mbox{\footnotesize\bf T}}
=D+\ell'$.
By construction, $D$ is the total space of the family of plane
quartic curves residual to $\ell$.
Finally, we have
$\L:=\O(D)=\O(5\cdot y-(y-x))=\O(4y+x)$.

Using \schub\cite{schub} we may compute
$
\int_{G_{2,\ell}}(c_6\eml{(2,2)}{\L}/2=
1,185$
(see the appendix) and find the number
$$\bf17\!,\!601\!,\!000\rm=tg_{6,5} - 609250 - 1185\times2875$$
of $irreducible$ plane rational quintic curves contained in a
generic  5\up{ic} threefold. The  1\up{st} correction is due to
$conic+cubic$ and the 2\up{nd} to $line+binodal$ $quartic$.

\section{Final comments}
An additional difficulty appears for the case of $7-$fold tangent
hyperplanes. Indeed, for a general $7-$dimensional linear
system, we'd expect $\Sigma(2^{[7]};S)$ to receive contributions
from $\Sigma(3(2),2;S),\ \Sigma(3,2^{[3]};S),\ \Sigma(3(2)';S)$,
so that a na\"{\i}ve count would predict
$$\mbox{\footnotesize$tg_7:=(\#\Sigma(2^{[7]};S)-
210\#\Sigma(3(2),2;S)-1260\#\Sigma(3,2^{[3]};S)/6)-
30\#\Sigma(3(2)';S))/7!,$}$$
\n\normalsize where $\Sigma(3(2)';S)$ denotes a cycle supported
on the set of $(C,y_1,\dots,y_7)$ such that $C\!\in\! S$ has a
triple point $y_1$ with the infinitely near double point $y_2$
presenting a branch tangent to the exceptional line over $y_1$.
However, barring some computational error, in fact the rhs did
not yield an integer for any of the examples we've experimented
with. This seems to indicate that $\Sigma(2^{[7]};S)$ may not be
reduced at some of the points involving singularities worse than
nodes. In fact, the argument of Prop.\ref{reduced} does not
apply. This would imply that the coefficients $210$, $1260$ and
$30$, postulated by the na\"{\i}ve count of permutations, must be
modified.

For $n\geq8$, we face the intrusion of a component of wrong
dimension in $\Sigma(2^{[n]};S)$ due to 4--fold points. In
this case, the technique of residual intersections
might shed some light.
\section{Appendix: computations}

\baselineskip6pt
\footnotesize
\begin{verbatim}
###CUT HERE FOR MAPLE
with(schubert):with(SF):

#PRINCIPAL PARTS of order n,
# f =cotg,d=linebundle
princ:= proc(n,f,d)local i:d&*sum('symm(i,f)',i=0..n):end:

whichmon:=proc(f,vars)local i,v,z:
   z:=expand(f):
      if type(z,`*`)or type(z,`^`)or type(z,`name`) then
	 v:=[seq(vars[i]=1,i=1..nops(vars))]:
	 RETURN(f/subs(v,f)):
      else ERROR(`invalid arg`)
      fi:end:

#SUBS EXACT MONOMIAL RELATIONS
submonpol:=proc(f,vars,rels)local z,i,j,term,mono,temp:
   z:=expand(f):temp:=0:
   if type(z,`+`)then
      for i to nops(z)do
	 term:=op(i,z):
	 mono:=whichmon(",vars):
	 for j to nops(rels) while mono<>lhs(rels[j])do od:
	 if j<=nops(rels)then
	    temp:=temp+term/mono*rhs(rels[j])
	 else temp:=temp+term:
	 fi:
      od:
   elif type(z,`*`)or type(z,`^`)or type(z,`name`)then
      term:=z:
      mono:=whichmon(term,vars):
      for j to nops(rels)while mono<>lhs(rels[j])do od:
      if j<=nops(rels)then
	 temp:=temp+term/mono*rhs(rels[j])
      else temp:=temp+term
      fi:
   fi:
RETURN(temp)end:

#KILL TERMS IN VARS OF TOTDEG>DIM
dimsimpl:= proc(x,vars,degs,dim)local i,j,temp,par,n:
  temp:=expand(x):
   if type(temp,`+`)then
      par:=0:n:=nops(temp):
      for i to n do
	 op(i,temp):
	 degree(collect(subs([seq(vars[j]=t_^degs[j]*vars[j],
	    j=1..nops(vars))],"),t_),t_):
	 if "<=dim then par:=par+"":fi:
      od:
      temp:=par:
   else
      degree(collect(subs([seq(vars[j]=t_^degs[j]*vars[j],
	 j=1..nops(vars))],"),t_),t_):
      if ">dim then temp:=0 fi:
   fi:
RETURN(temp):end:

simplification:=proc () local i, j, n, z, zz:
   n:= args[1]: z:=args[2]:
   if nargs=3 and type(args[3],set) then
      zz:=args[3] else zz:={n}
   fi:
   for i  from n by -1 to 2 do for j to n+1-i do
      if 2 < degree(collect(z,e[j,i]),e[j,i]) then
	 z:=rem(collect(z,e[j,i]),relexc.i.j,e[j,i]):
	 zz:=zz union{i+j-1}:
      fi:
   od:od:
   if opt_=5 then
      for i in zz do
	 if 2 < degree(collect(z,y.i),y.i) then
	    z:=rem(z,rely.i,y.i):
	 fi:
      od:
   else
      for i in zz do z:=dimsimpl(z,var0.(i),deg1,2):od:
   fi :
RETURN(z)end:

#MAIN PROCEDURE FOR PUSHFORWARD {n}->{n-1}
push:= proc(n,f)
  local z,z0,z2,zz,mons,i,j,i1,j1,temp,varn,var0,degn,dd:
  option remember:
     if opt_=5 and type(relpush_5,set)=false then relpush_5:={}:fi:
     if opt_<>5 and type(relpush_,set)=false then relpush_:={}:fi:
     if n=1 then
	if opt_=5 then
	   rem(f,rely1,y1):z:=coeff(collect(",y1),y1,2):
	else
	   subs([seq(var1[i]=0,i=1..nops(var1))],f):
	   f-":dimsimpl(",var0.n,deg1,2):
	   submonpol(",var1,{c[1,2]=chi,c[1,1]^2=k2,h[1]^2=d,
	      h[1]*c[1,1]=hk}):
	   z:=submonpol(",var1,{c[1,1]=0,h[1]=0})
	fi:
     RETURN(z):
     else
	convert(var.(n-1),set) minus convert(var.(n-2),set):
	var0:=[op(")]:
	degn:=[seq(1,i=1..nops(var0))]:
	if opt_<>5 then
	  var0:=[op(var0),c[n-1,2]]:degn:=[op(degn),2]
	fi:
	varn:=[seq(p.n.2 &^* var0[i],i=1..nops(var0))]:
	subs([e[1,n]=0,seq(varn[i]=0,i=1..nops(varn))],f):
	z:=collect(f-",e[1,n]):
	if 2<degree(z,e[1,n]) then z:=rem(z,relexc.n.1,e[1,n]):
	fi:
	z:=collect(z,e[1,n]):
	z:=collect(z-e[1,n]*coeff(z,e[1,n],1),e[1,n]):
	if z<>0 then
	   z0:=coeff(z,e[1,n],0):
	   if z0<>0 then
	      simplification(n,z0):z0:=collect(",e[2,n-1]):
	      z0:=z0-e[2,n-1]*coeff(z0,e[2,n-1],1)
	   fi:
	   zz:=collect(subs([seq(varn[i]=t_*varn[i],
	      i=1..nops(varn))],z0),t_):
	   dd:=degree(zz,t_):
	   temp:=0:
	   for i from dd by -1 to 1 do
	      z0:=expand(coeff(zz,t_,i)):
	      if z0 <> 0 then
		 if (opt_=5 and type(relpush_5.n.i,list)=false
		    or (opt_<>5 and type(relpush.n.i,list)=false)) then
		    if opt_=5 then
		       print(`BUILD RELPUSH_5`.n.i):
		    elif opt_<>5 then
		       print(`BUILD RELPUSH`.n.i):
		    fi:
		    mons :=monomials(i,var0,degn):
		    z2:={}:
		    for j to nops(mons) do
		       dimsimpl(mons[j],var0.(n-1),deg1,2):
		       if degree(collect(",e[1,n-1]),e[1,n-1])<>1 and
			  member(true,{seq(type("/e[j1,n-j1]^3,
			     polynom),j1= 1..n-2)})=false then
			  z2:=z2 union {mons[j]}
		       fi:
		    od:
		    mons:=[seq(p.n.2&^* z2[j]=push(n-1,z2[j]),j=
		       1..nops(z2))]:
		    if opt_=5 then relpush_5.n.i:=mons:
		    else relpush.n.i:=mons:
		    fi:
		 elif opt_=5 and member([n,i],relpush_5)=false then
		    print(`USING RELPUSH_5`.n.i.` BUILT BEFORE`):
		    relpush_5:=relpush_5 union{[n,i]}:
		 elif opt_<>5 and member([n,i],relpush_)=false then
		    print(`USING RELPUSH`.n.i.` BUILT BEFORE`):
		    relpush_:=relpush_ union{[n,i]}
		 fi:
		 if opt_=5 then mons:=relpush_5.n.i:
		 else mons:=relpush.n.i:
		 fi:
		 z0:=submonpol(z0,varn,mons)
	      fi:
	      temp:=temp+z0:
	   od:
	   z0:=temp:
	   z2:=-coeff(z,e[1,n],2):
	   for i to nops(varn) while z2 <> 0 do
	      z2:=collect(z2,varn[i]):
	      if degree(z2,varn[i]) <> 0 then
		 z2:=rem(z2,varn[i]-var0[i],varn[i]):
	      fi:
	   od:
	   if z2<>0 then z2:=simplification(n-1,z2) fi:
	   z:=z0+z2:
	fi:
	RETURN(z)
     fi:
end:  #of push

# CALCULATIONS
for opt_ in[5,0]do
   if opt_=5 then
      grass(3,5,x,all):
      Grass(g,1,Qx,y,all):
      omega1:=dual(g[tangentbundle_]):
      rely1:=chern(3,"):var1:=[y1]:deg1:=[1]:
      variety(S1,dim=8,vars=var1,degs=deg1):
   else
      var1:=[c[1,1],h[1],c[1,2]]:deg1:=[1,1,2]:
      variety(S1,dim=2,vars=var1,degs=deg1):
   fi:
   var01:=var1:
   for n to 6 do
      if n=1 then
	 if opt_=5 then
	    DIM:=3:         #ONLY FOR THE SAKE OF RANKS...
	    L:=o(4*y1+x1):  #FOR BINODAL 4ICS
	    princ(1,omega1,L):chern(3,"):FB1:=rem(",rely1,y1):
	    print(`done FB1`):
	    DIM:=6:M:=o(m*y1):DIM:=3:
	 else
	    opt_:=0:goto(S1):bundle(2,c):
	    omega1:=subs([c1=c[1,1],c2=c[1,2]],"):   #COTANGENT BUNDLE
	    M:=o(h[1]):
	 fi:
	 princ(1,omega1,M):
	 if opt_=5 then chern(3,"):else chern("):
	 fi:
	 F1:=simplification(n,"):print(`done F1`):
	 if opt_=5 then
	    DIM:=6:         #ONLY FOR THE SAKE OF RANKS...
	 fi:
	 princ(2,bundle(2,c),M):
	 subs([seq(c.i=chern(i,omega1),i=1..2)],"):
	 if opt_=5 then chern(6,"):else chern("):
	 fi:
	 E_31:=simplification(1,"):print(`done E_31`):
      elif n>=2 then
	 var0.n:=[y.n]:
	 var.n:=[seq(y.j,j=1..n),
	    seq(seq(e[j,k],j=1..n-k+1),k=2..n-1),e[1,n]]:
	 deg.n:=[seq(1,j=1..n),
	    seq(seq(1,j=1..n-k+1),k=2..n-1),1]:
	 if opt_<>5 then
	    var0.n:=[c[n,1],h[n],c[n,2]]:
	    var.n:=[seq(c[j,1],j=1..n),
	       seq(c[j,2],j=1..n),op(subs([seq(y.j=h[j],j=1..n)],var.n))]:
	    deg.n:=[seq(1,j=1..n),seq(2,j=1..n),op(deg.n)]
	 fi:
	 if opt_=5 then rely.n:=subs(y1=y.n,rely1): fi:
	 variety(S.n,dim=6/5*opt_+2*n,vars=var.n,degs=deg.n):
	 morphism(p.(n).2,S.n,S.(n-1),subs([seq(var0.(n-1)[k]=var0.(n)[k],
	    k=1..nops(var.01)),seq(e[n-k,k]=e[n-k+1,k],k=2..n-1)],
	    var.(n-1))):print(`built S`.n):
	 DIM:=3:    #OK since ranks<=3
	 omega.n:=((p.n.2)&^*(omega.(n-1)))&*o(e[1,n])+o(-e[1,n])-1:
	 chern(3,omega.n):       #Will set=0 since rk.omega=2
	 print(`DONE OMEGA`.n):
	 relexc.n.1:=rem(",rely.n,y.n):
	 for i from n-1 by -1 to 2 do
	    relexc.i.(n+1-i):=(p.n.2)&^*(relexc.i.(n-i)):
	    if degree(collect(relexc.n.1,e[n+1-i,i]),e[n+1-i,i])>2 then
	       relexc.n.1:=rem(relexc.n.1,",e[n+1-i,i])
	    fi:
	 od:
	 M:=collect((p.n.2)&^*M&*o(-2*e[1,n]),t):        #Adjust M
	 princ(1,omega.n,M):
	 if opt_=5 then chern(3,"):else chern("): fi:
	 F.n:=simplification(n,"):print(`done F`.n):
	 if n=2 then #E_32
	    collect(M&*o(-e[1,2]),t):princ(1,omega2,"):
	    if opt_=5 then
	       E_32:=chern(3,"):
	       L:=collect((p.n.2)&^*L&*o(-2*e[1,n]),t):        #Adjust L
	       princ(1,omega.n,L):chern(3,"):
	       FB.n:=simplification(n,"):print(`done FB`.n):
	    else E_32:=chern("):
	    fi:
	    E_32:=simplification(n,E_32):print(`done E_32`):
	    E_3_2:=rem(collect("*e[1,2],e[1,2]),relexc21,e[1,2]):
	    print(`done E_3_2`):
	 fi:#n=2
	 if n=3 then #E_33
	    collect(M&*o(-e[2,2]),t):princ(1,omega3,"):
	    if opt_=5 then chern(3,"):else chern("): fi:
	    E_33:=simplification(n,"):print(`done E_33`):
	 fi: #E_33
      fi:
   od:
   if opt_=5 then
      1:
      for i from 6 by -1 to 1 do F.i:=push(i,"*F.i):print(i)
      od:
      E_33:=push(3,E_33):E_32:=push(2,E_32*E_33):
      E_3_2:=push(2,E_3_2):e_322:=push(1,E_31*E_32):
      e_3_2:=push(1,E_31*E_3_2):
      tg.6:=integral(Gx,F1-30*e_3_2-90 *e_322)/6!:
      lprint(`#4-coplanar lines in 4ic 3fld: `,subs(m=4,tg.6)):
      lprint(`#6-nodal plane sections of 5ic 3fold: `,subs(m=5,tg.6)):
      lprint(`#binodal plane 4ic residul to line in 5ic 3fld: `,
	1/2*integral(Gx,(x1^2-x2)^2*push(1,FB1*push(2,FB2)))):
      lprint(`#6-nodal IRREDUCIBLE plane sections of 5ic 3fold: `,
	subs(m=5,tg.6) - 609250 - 1185*2875):
      #4-COPLANAR LINES VIA FANO
      DIM:=3:
      for i to 3 do relm.i.1:=chern(3,5-bundle(2,z.i)-o(m.i.1)):
	 for j to 2 do
	   z.i.j:=chern(j,bundle(3,x)-o(m.i.1)):
       od od:
       chern(3,symm(2,bundle(2,z3))&*o(2*x1-z11-z21)):
       rem(",relm31,m31):
       I_3:=coeff(",m31,2):
       DIM:=4:
       chern(4,symm(3,bundle(2,z2))&*o(x1-z11))*I_3:
       rem(",relm21,m21):
       I_2:=coeff(",m21,2):
       DIM:=5:
       I_1:=chern(5,symm(4,bundle(2,z1)))*I_2:
       I_1:=rem(I_1,relm11,m11):
       I_1:=coeff(I_1,m11,2):
       integral(Gx,")/4!:
    else
       for j to 6 do for i from 0 to 3 do
	  F.j.i:=coeff(collect(F.j,t),t,i):
	  if j<=3 then E_3.j.i:=coeff(collect(E_3.j,t),t,i):fi:
	  if j=2 then E_3_2.(i):=coeff(collect(E_3.2,t),t,i)*e[1,2]:fi:
       od:od:
       for i from 1 to 3 do F6.i:=push(6,F6.i):od:
       ftg6:=ftg5*F62:
       for j5 from 1 to 3 do print(`j5=`.j5):
	  a5:=push(5,`F`.5.j5*F63):
	  #dim 10-j5-1=9-j5<=8 ok

	  for j4 from 3-j5 to 3 do
	     a4:=push(4,`F`.4.j4*a5):
	     # 6>=dim 9-j5-j4 >=0

	      for j3 from max(0,5-j5-j4) to 3 do
		 a3:=push(3,`F`.3.j3*a4):
		 #4>=dim 9-j5-j4-j3 >=0

		 for j2 from max(0,7-j5-j4-j3) to min(3,9-j5-j4-j3) do
		    a2:=push(2,`F`.2.j2*a3):
		    #2>=dim 9-j5-j4-j3-j2
		    j1 :=9-j5-j4-j3-j2 :
	  lprint(`j5=`.j5,` j4=`.j4,` j3=`.j3,` j2=`.j2,` j1=`.j1):
		    ftg6:=ftg6+push(1,`F`.1.j1*a2):
	od:od:od:od:
###########
  for i from 1 to 3 do F5.i:=push(5,F5.i):od:
  ftg5:=ftg4*F52:
  for j4 from 1 to 3 do
     a4:=push(4,`F`.4.j4*F53):
     #dim 8-j4-1=7-j4<=6 ok

     for j3 from 3-j4 to 3 do
	a3:=push(3,`F`.3.j3*a4):
	#4>=dim 7-j4-j3 >=0

	 for j2 from max(0,5-j4-j3 ) to min(3,7-j4-j3) do
	    a2:=push(2,`F`.2.j2*a3):
	    #2>=dim 7-j4-j3-j2>=0
	    j1 :=7-j4-j3-j2 :
	lprint(`j4=`.j4,` j3=`.j3,` j2=`.j2,` j1=`.j1):
	    ftg5:=ftg5+push(1,`F`.1.j1*a2):
  od:od:od:
###########
  for i from 1 to 3 do F4.i:=push(4,F4.i):od:
  ftg4:=ftg3*F42:
  for j3 from 1 to 3 do
     a3:=push(3,`F`.3.j3*F43):
     #4>=dim 5-j3 >=0

     for j2 from max(0,3-j3) to min(3,5-j3) do
	a2:=push(2,`F`.2.j2*a3):
	#2>=dim 5-j3-j2>=0
	j1 :=5-j3-j2 :
       lprint(`j3=`.j3,` j2=`.j2,` j1=`.j1):
	ftg4:=ftg4+push(1,`F`.1.j1*a2):
  od:od:
##########
  for i from 1 to 3 do F3.i:=push(3,F3.i):E_33.i:=push(3,E_33.i):
  od:
  ftg3:=ftg2*F32:
  for j2 from 1 to 3 do
    a2:=push(2,`F`.2.j2*F33):
    #2>=dim 3-j2 >=0
    j1 :=3-j2 :
    lprint(` j2=`.j2,` j1=`.j1):
    ftg3:=ftg3+push(1,`F`.1.j1*a2):
  od:
  e_322:=e_32*E_332:
  for j2 from 1 to 3 do
    a2:=push(2,`E_3`.2.j2*E_333):
    #2>=dim 3-j2 >=0
    j1 :=3-j2:
    lprint(`j2=`.j2,` j1=`.j1):
    e_322:=e_322+push(1,`E_3`.1.j1*a2):
  od:
##########
  for i from 1 to 3 do F2.i:=push(2,F2.i):
    E_32.i:=push(2,E_32.i):
    E_3_2.(i):=push(2,E_3_2.(i)):
  od:
  ftg2:=ftg1*F22+push(1,`F`.1.1*F23):
  e_32:=e_3*E_322+push(1,`E_3`.1.1*E_323):
  e_3_2:=e_3*E_3_21+push(1,`E_3`.1.1*E_3_22)+
    push(1,`E_3`.1.0*E_3_23):
  ftg1:=push(1,F12):e_3:=push(1,E_312):
  for n to 3 do tg.n:=ftg.n /n! :od:
  tg.4:=(ftg4-6*e_3)/4!:
  tg.5:=(ftg5-30*e_32)/5!:
  tg.6:=(ftg6-30*e_3_2-90*e_322)/6!:
  p2:=proc(m,tg)subs([chi=3,d=m^2,hk=-3*m,k2=9],  tg):end:
  p1xp1:=proc(m1,m2,tg)subs([chi=4,k2=8,hk=-2*m2-2*m1,
    d=m1*m2*2],tg):end:
  K_3:=proc(g)subs([chi=24,k2=0,hk=0,d=2*g-2],`tg`.g):end:
  for i from 3 to 6 do
    lprint(`# `.i.`-nodal hyperplane sections of K3-sfce in P`.
      i.`: `,K_3(i)):
  od:
  for i from 4 to 6 do lprint(`# `.i.`-nodal plane quartics through `.
     (14-i).` general points: `,p2(4,tg.i)):
   od:
fi:
od:
#CUT HERE FOR MAPLE \end{verbatim}
\vfill\eject

{\baselineskip5pt\parskip0pt
\begin{obeylines}
\n Departamento de Matem\'atica
\n Universidade Federal de Pernambuco
\n Cidade Universit\'aria
\n 50670--901 Recife--Pe--Brazil
\n email: israel@dmat.ufpe.br
\end{obeylines}
}
\enddocument